\documentclass{article}%
\usepackage{amsmath}
\usepackage{amsfonts}
\usepackage{amssymb}
\usepackage{graphicx}%
\providecommand{\U}[1]{\protect\rule{.1in}{.1in}}
\newtheorem{theorem}{Theorem}

\newtheorem{corollary}[theorem]{Corollary}

\newtheorem{definition}[theorem]{Definition}

\newtheorem{lemma}[theorem]{Lemma}

\newtheorem{proposition}[theorem]{Proposition}
\newtheorem{remark}[theorem]{Remark}

\newenvironment{proof}[1][Proof]{\noindent\textbf{#1.} }{\ \rule{0.5em}{0.5em}}
\begin{document}

\title{Characterization of the attractor for nonautonomous reaction-diffusion
equations with discontinuous nonlinearity}
\author{Jos\'{e} Valero\\{\small Universidad Miguel Hernandez de Elche,}\\{\small Centro de Investigaci\'{o}n Operativa,}\\{\small Avda. Universidad s/n,}\\{\small 03202-Elche (Alicante), Spain }\\{\small jvalero@umh.es }}
\date{}
\maketitle

\begin{abstract}
In this paper, we study the asymptotic behavior of the solutions of a
nonautonomous differential inclusion modeling a reaction-diffusion equation
with a discontinuous nonlinearity.

We obtain first several properties concerning the uniqueness and regularity of
non-negative solutions. Then we study the structure of the pullback attractor
in the positive cone, showing that it consists of the zero solution, the
unique positive nonautonomous equilibrium and the heteroclinic connections
between them, which can be expressed in terms of the solutions of an
associated linear problem.

Finally, we analyze the relationship of the pullback attractor with the
uniform, the cocycle and the skew product semiflow attractors.

\end{abstract}

\bigskip

\textbf{Keywords: }differential inclusions, reaction-diffusion equations,
pullback attractors, nonautonomous dynamical systems, multivalued dynamical
systems, structure of the attractor

\textbf{AMS Subject Classification (2010): }35B40, 35B41, 35B51, 35K55, 35K57

\section{Introduction}

The study of the structure of the global attractor for dynamical systems is an
important question as it leads to a deep understanding of the dynamics of the
solutions of differential equations. In the single-valued autonomous case,
there are good examples in which the dynamics inside the attractor has been
fully described (see for example \cite{BrunoswkyFiedler88},
\cite{BrunoswkyFiedler89}, \cite{Fiedler96}, \cite{Hale}, \cite{Henry85},
\cite{Rocha88}, \cite{Rocha91}). In the multivalued autonomous case, it is
more difficult to carry out the study of the structure of attractors.
Nevertheless, several results have been published in this direction (see
\cite{ARV06}, \cite{CabCarvMarVal}, \cite{KKV14}, \cite{KKV15}). In all these
papers, the dynamical system is of gradient type, so the attractor is
described by means of the set of stationary points and its unstable manifold.
Moreover, in some of them, the dynamics is described in detail in terms of the
stationary points and their heteroclinic connections. The Chaffee-Infante
equation \cite{Henry85} is a paradigmatic example in which the dynamics has
been completely understood.

In the nonautonomous case, the problem is more complex as stationary solutions
do not exist in general in a classical sense. Instead, we need to replace them
by a special type of bounded complete trajectories, which play the role of
\textquotedblleft nonautonomous equilibria\textquotedblright. In this
direction, the existence of a complete bounded positive non-degenerate
trajectory is a key fact which has been proved successfully for parabolic
equations of certain types (see \cite{CLR}, \cite{LRBS}, \cite{LangaSuarez},
\cite{RRV}, \cite{RBVidal}). Moreover, in the special case of the
nonautonomous Chaffee-Infante equation it is proved that there exist a finite
number of nonautonomous equilibria, which are analogous to the equilibria in
the autonomous situation, and that the omega and alpha limits of any bounded
complete trajectory belong to the set of these equilibria when we consider all
the equations generated by the translations of the nonautonomous term
\cite{BCV19}. Therefore, the pullback attractor is characterized in terms of
the nonautonomus equilibria and their heteroclinic connections, leading to a
gradient structure of the attractor of the associated skew product semiflow.

In this paper, we study the dynamical system generated by the differential
inclusion%
\begin{equation}
\left\{
\begin{array}
[c]{l}%
\dfrac{\partial u}{\partial t}-\dfrac{\partial^{2}u}{\partial x^{2}}\in
b(t)H_{0}(u)+\omega(t)u,\text{ on }(\tau,+\infty)\times(0,1),\\
u(t,0)=u(t,1)=0,\\
u(\tau,x)=u_{\tau}(x),\text{ }x\in\left(  0,1\right)  ,
\end{array}
\right.  \label{Eq}%
\end{equation}
where
\begin{equation}
H_{0}(u)=\left\{
\begin{array}
[c]{ll}%
-1, & \text{if }u<0,\\
\left[  -1,1\right]  , & \text{if }u=0,\\
1, & \text{if }u>0,
\end{array}
\right.  \label{Heaviside}%
\end{equation}
is the Heaviside function, so we have in fact a differential equation with a
discontinuous nonlinearity. This dynamical system is both nonautonomous and
multivalued. Well known models like combustion in porous media \cite{FeNo},
the conduction of electrical impulses in nerve axons (see \cite{Terman},
\cite{Terman2}) or the surface temperature on Earth (see \cite{Budyko},
\cite{DiazHernandezTello}) have discontinuities of this type. This inclusion
is also important because it is the limit of a sequence of Chaffee-Infante
problems which has undergone all the bifurcation cascade of this type of
problems (see \cite{ARV06} for more details), so it is natural to expect that
it should inherit the structure of the attractor. In fact, in the autonomous
case (that is, when the functions $b\left(  \text{\textperiodcentered}\right)
,\omega\left(  \text{\textperiodcentered}\right)  $ are constants), it is
proved in \cite{ARV06} that (\ref{Eq}) has an infinite, but countable, number
of equilibria (each of which is related to a corresponding one in the
Chaffee-Infante equation) and that the attractor is described by them and
their heteroclinic connections. Moreover, when $\omega$ is equal to $0$ it is
shown that some of the connections of the Chaffee-Infante equation are also
true for problem (\ref{Eq}). The hypothesis is that the connections are the
same and then the structure of the attractors coincide. However, this still
remains an open problem. We observe also that for a related inclusion modeling
the climate on Earth several results on bifurcations of steady states were
proved in \cite{BensidDiaz}, \cite{BensidDiaz2}.

In this paper, we continue the study of the structure of the pullback
attractor for problem (\ref{Eq}) initiated in \cite{CVL19}, where the
existence of two special non-degenerate bounded complete trajectories
containing the pullback attractor was established (the positive and negative
nonautonomous equilibria). More precisely, we give a complete characterization
of the pullback attractor in the positive cone, that is, when we consider only
non-negative solutions.

First, using the maximum principle as the main tool we prove the following
important facts:

\begin{itemize}
\item For any non-negative initial condition not equal to zero there exists a
unique non-negative solution, which becomes positive instantaneously.

\item If the initial condition is equal to zero, then the unique non-negative
solution backwards in time is the zero solution.

\item If the initial condition is non-negative and not equal to zero, belongs
to the space $H_{0}^{1}\left(  0,1\right)  $ and is not positive, then there
cannot exist a non-negative solution backwards in time.

\item If the initial condition is non-negative and belongs to $L^{2}\left(
0,1\right)  \backslash H_{0}^{1}\left(  0,1\right)  $, then there cannot exist
a non-negative solution backwards in time.
\end{itemize}

Second, using these results we obtain that the pullback attractor in the
positive cone consists of two nonautonomous equilibria (the zero solution and
the positive nonautonomous equilibrium denoted by $\xi_{M}$) and all the
bounded complete trajectories which connect them heteroclinically. Moreover,
these connections always go from $0$ to $\xi_{M}$ and have the form%
\[
\varphi(t)=\left\{
\begin{array}
[c]{c}%
0\text{ if }t\leq\tau,\\
u\left(  t\right)  \text{ if }t\geq\tau,
\end{array}
\right.
\]
where $u\left(  \text{\textperiodcentered}\right)  $ is the unique solution to
the associated linear problem of (\ref{Eq}) in which the right-hand side is
equal to $b(t)+\omega(t)u\left(  t\right)  $ and we choose the initial
condition $u_{\tau}=0$. It is important to point out that these results are
also new in the autonomous case. We observe also that the structure of the
pullback attractor in the positive cone is the same as in the associated
ordinary differential inclusion in which the diffusion term $-\dfrac
{\partial^{2}u}{\partial x^{2}}$ is replaced by the linear function $\lambda
u$, $\lambda>0$ \cite{CVL16}.

Additionally, we obtain the following properties of the positive nonautonomous
equilibrium $\xi_{M}$, which were left as an open problem in \cite[Remark
2]{CVL19}:

\begin{itemize}
\item The solution starting at any point $\xi_{M}\left(  s\right)  $,
$s\in\mathbb{R}$, is unique in the class of all non-negative solutions. As a
consequence, $\xi_{M}$ is the unique non-negative non-degenerate bounded
complete trajectory at $-\infty$.

\item If $b,\omega$ are more regular ($b,\omega\in W_{loc}^{1,2}\left(
\mathbb{R}\right)  $), then the solution starting at any point $\xi_{M}\left(
s\right)  $, $s\in\mathbb{R}$, is unique in the class of all solutions. Thus,
$\xi_{M}$ is the unique non-degenerate bounded complete trajectory at
$-\infty$.
\end{itemize}

Moreover, increasing the regularity of the functions $b,\omega$ we increase
the regularity (in both time and space) of non-negative solutions as well. In
this way, we are able to prove that the attractor in the positive cone is as
regular as we desire.

Third, we obtain the existence of other types of nonautonomous attractors such
as the cocycle attractor, the uniform attractor and the skew product semiflow
attractor and describe their relationship. In particular, we obtain that the
associated skew product semiflow possesses in the positive cone a global
attractor with a gradient structure.

\section{Setting of the problem}

We consider the problem
\begin{equation}
\left\{
\begin{array}
[c]{l}%
\dfrac{\partial u}{\partial t}-\dfrac{\partial^{2}u}{\partial x^{2}}\in
b(t)H_{0}(u)+\omega(t)u,\text{ on }(\tau,\infty)\times(0,1),\\
u(t,0)=u(t,1)=0,\\
u(\tau,x)=u_{\tau}(x),\text{ }x\in\left(  0,1\right)  ,
\end{array}
\right.  \label{Incl}%
\end{equation}
where $b:\mathbb{R}\rightarrow\mathbb{R}^{+},$ $\omega:\mathbb{R}%
\rightarrow\mathbb{R}^{+}$ are continuous functions such that
\begin{equation}
0<b_{0}\leq b\left(  t\right)  \leq b_{1},\ 0\leq\omega_{0}\leq\omega\left(
t\right)  \leq\omega_{1}, \label{Cond1}%
\end{equation}
and $H_{0}$ is the Heaviside function given in (\ref{Heaviside}).

Let $H=L^{2}(0,1)$, $V=H_{0}^{1}(0,1)$. The norm in $H$ will be denoted by
$\left\Vert \text{\textperiodcentered}\right\Vert $, whereas the norm in any
other Banach space $X$ will be denoted by $\left\Vert
\text{\textperiodcentered}\right\Vert _{X}$. $L_{loc}^{p}(\tau,+\infty;X)$,
$p\geq1$ (respectively, $W_{loc}^{k,p}(\tau,+\infty;X)$, $p\geq1,\ k\in
\mathbb{N}$) will stand for the space of functions $g$ such that $g\in
L^{p}(\tau,T;X)$ (respectively, $W^{k,p}(\tau,T;X)$) for any $T>\tau.$

For a metric space $Y$ denote by $dist_{Y}(A,B)=\sup_{a\in A}\inf_{b\in B}%
\rho_{Y}(a,b)$ the Hausdorff semidistance from the set $A$ to the set $B$,
where $\rho_{Y}$ is the metric in $Y.$

Let $A:D(A)\rightarrow H,\ D(A)=H^{2}(0,1)\cap V,$ be the operator
$A=-\dfrac{d^{2}}{dx^{2}}$ with Dirichlet boundary conditions. This operator
is the generator of a $C_{0}$-semigroup $T(t)=e^{-At}$. Also, we define the
multivalued operator $R:\mathbb{R}\times H\rightarrow2^{H}$ (with $2^{H}$
being the set of all subsets of $H$) by%
\[
R(t,u)=\left\{  y\in H:\text{ }y\left(  x\right)  \in b(t)H_{0}\left(
u\left(  x\right)  \right)  +\omega(t)u\left(  x\right)  \text{, a.e. on
}\left(  0,1\right)  \right\}  .
\]
It is known \cite[Lemma 2]{CVL19} that this operator has nonempty, closed,
bounded and convex values.

We rewrite (\ref{Incl}) in the abstract form%
\[
\left\{
\begin{array}
[c]{c}%
\dfrac{du}{dt}+Au\in R(t,u\left(  t\right)  ),\ t\in\left(  \tau
,+\infty\right)  ,\\
u\left(  \tau\right)  =u_{\tau}.
\end{array}
\right.
\]

\begin{definition}
(Strong solution) Let $u_{\tau}\in H$. The function $u\in C([\tau,+\infty),H)$
is said to be a strong solution to problem (\ref{Incl}) if:

\begin{enumerate}
\item $u(\tau)=u_{\tau};$

\item $u\left(  \text{\textperiodcentered}\right)  $ is absolutely continuous
on $[T_{1},T_{2}]$ for any $\tau<T_{1}<T_{2}$ and $u\left(  t\right)  \in
D(A)$ for a.a. $t\in\left(  T_{1},T_{2}\right)  ;$

\item There exists a function $r\in L_{loc}^{2}(\tau,+\infty;H)$ such that
$r\left(  t\right)  \in R(t,u\left(  t\right)  )$ for a.a. $t\in\left(
\tau,+\infty\right)  $ and%
\begin{equation}
\frac{du}{dt}+Au(t)=r(t)\text{ for a.a. }t\in\left(  \tau,+\infty\right)  ,
\label{Equality}%
\end{equation}
where the equality is understood in the sense of the space $H.$
\end{enumerate}
\end{definition}

\begin{definition}
(Mild solution) Let $u_{\tau}\in H$. The function $u\in C([\tau,+\infty),H)$
is said to be a mild solution to problem (\ref{Incl}) if there exists a
function $r\in L_{loc}^{2}(\tau,+\infty;H)$ such that $r\left(  t\right)  \in
R(t,u\left(  t\right)  )$ for a.a. $t\in\left(  \tau,+\infty\right)  $ and%
\[
u(t)=e^{-A(t-\tau)}u_{\tau}+\int_{\tau}^{t}e^{-A(t-s)}r(s)ds\text{ for }%
\tau\leq t<+\infty.
\]

\end{definition}

It is known \cite[Theorem 1]{CVL19} that for any $u_{\tau}\in H$ there exists
at least one strong solution to problem (\ref{Incl}). For any $u_{\tau}\in H$
the function $u\in C([\tau,+\infty),H)$ is a strong solution to problem
(\ref{Incl}) if and only if it is a mild solution \cite[Corollary 4]{CVL19}.
Therefore, from now on we will speak just about solutions of (\ref{Incl}).

The main aim of this paper is to establish a\ precise characterization of the
pullback attractor for non-negative solutions to problem (\ref{Incl}) in terms
of the bounded complete trajectories connecting the stationary solutions.

\section{About positive solutions}

An element $v\in H$ is said to be non-negative (denoted by $v\geq0$) if
$v\left(  x\right)  \geq0$ for a.a. $x\in\left(  0,1\right)  $. Let us
consider non-negative solutions to problem (\ref{Incl}). Under assumption
(\ref{Cond1}) it is known \cite[Corollary 5]{CVL19} that for any $u_{\tau}\in
H$ satisfying $u_{\tau}\geq0$ there exists at least one strong solution
$u\left(  \text{\textperiodcentered}\right)  $ such that $u\left(  t\right)
\geq0$ for any $t\geq\tau$.

An element $v\in V$ is said to be positive if $v(x)>0$ for all $x\in\left(
0,1\right)  $. We shall show in this section that for initial conditions
$u_{\tau}$ that are non-negative and not identically equal to zero any
solution $u$ to problem (\ref{Incl}) satisfies that $u(t)$ is positive for any
$t>0.$ Also, for such initial conditions the solution is unique among
non-negative solutions. Moreover, we establish that if either the initial
condition is not positive and not equal to zero or it belongs to $H\backslash
V$, then the solution does not exist backwards in time.

We consider the following auxiliary problem:%
\begin{equation}
\left\{
\begin{array}
[c]{l}%
\dfrac{\partial u}{\partial t}-\dfrac{\partial^{2}u}{\partial x^{2}%
}=f(t),\text{ on }(\tau,\infty)\times(0,1),\\
u(t,0)=u(t,1)=0,\\
u(\tau,x)=u_{\tau}(x),\ x\in\left(  0,1\right)  ,
\end{array}
\right.  \label{Aux}%
\end{equation}
where $f\in L_{loc}^{1}(\tau,T;H)$.

For $u_{\tau}\in H$ the function \ $u\in C([\tau,+\infty),H)$ is called a
strong solution to problem (\ref{Aux}) if $u\left(  \tau\right)  =u_{\tau}$,
$u$ is absolutely continuous on any compact interval of $\left(  \tau
,+\infty\right)  $, $u\left(  t\right)  \in D(A)$ for a.a. $t\in\left(
\tau,+\infty\right)  $ and
\[
\frac{du}{dt}+Au(t)=f(t)\text{ in }H\text{ for a.a. }t\in\left(  \tau
,+\infty\right)  .
\]

The following lemma follows from a general result for equations governed by
subdifferential maps given for example in \cite[Theorem 3.6]{Brezis73} or
\cite[p.189]{Barbu}.

\begin{lemma}
\label{RegularityA}For any $f(\cdot)\in L_{loc}^{2}(\tau,+\infty;H)$,
$u_{\tau}\in H$, there exists a unique strong solution to problem (\ref{Aux})
satisfying
\begin{equation}
\sqrt{t}\frac{du}{dt}\in L_{loc}^{2}(\tau,+\infty;H),\,\,u\in L_{loc}^{2}%
(\tau,+\infty;V). \label{PropSolSelection}%
\end{equation}
Also, the map $t\mapsto\left\Vert u(t)\right\Vert _{V}$ is absolutely
continuous on any compact interval of $(\tau,+\infty)$.

If, moreover, $u_{\tau}\in V,$ then $\dfrac{du}{dt}\in L_{loc}^{2}\left(
\tau,+\infty;H\right)  $ and $t\mapsto\left\Vert u(t)\right\Vert _{V}$ is
absolutely continuous on any compact interval of $[\tau,+\infty)$.
\end{lemma}

\begin{remark}
\label{RegularityB}We observe that $u\in C([T_{1},T_{2}],H)\cap L^{\infty
}\left(  T_{1},T_{2};V\right)  $ implies that $u\in C([T_{1},T_{2}],V_{w})$,
where $V_{w}$ stands for the weak continuity in $V$ \cite[Lemma 1.4,
p.263]{Temam}. Then, as weak continuity in $V$ together with continuity in
norm gives strong continuity in the Hilbert space $V$, $u\in C([T_{1}%
,T_{2}],H)$ and $\left\Vert u(t)\right\Vert _{V}\in C([T_{1},T_{2}])$ imply
that $u\in C([T_{1},T_{2}],V)$. Hence, if $u_{\tau}\in H$ (respectively,
$u_{\tau}\in V$), then $u\in C(\left(  \tau,+\infty\right)  ,V)$
(respectively, $u\in C([\tau,+\infty),V)$).
\end{remark}

For $u_{\tau}\in H$ the function $u\in C([\tau,+\infty),H)$ is called a mild
solution to problem (\ref{Aux}) if
\[
u(t)=e^{-A(t-\tau)}u_{\tau}+\int_{\tau}^{t}e^{-A(t-s)}f(s)ds\text{ for }%
\tau\leq t<+\infty.
\]
By \cite[Lemma 4]{CVL19} the unique strong solution given in Lemma
\ref{RegularityA} is the unique mild solution to problem (\ref{Aux}).

We observe that by Remark \ref{RegularityB}\ any solution $u$ to problem
(\ref{Incl}) satisfies $u\in C((\tau,+\infty),V),$ so $V\subset C([0,1])$
gives $u\in C((\tau,+\infty),C([0,1]))$. Therefore, the solutions are
continuous in $\left(  \tau,+\infty\right)  \times\lbrack0,1]$. If we take an
initial condition in $V$, then $u\in C([\tau,+\infty),V),$ so the solutions
are continuous in $[\tau,+\infty)\times\lbrack0,1].$ Also, the derivative
$\dfrac{du}{dt}$ exists in the classical sense for a.a. $t>\tau$.

\begin{lemma}
Let $u\left(  \text{\textperiodcentered}\right)  $ be a non-negative solution
to problem (\ref{Incl}) with initial condition $u_{\tau}\in H$. Then%
\begin{equation}
\frac{\partial u}{\partial t}-\frac{\partial^{2}u}{\partial x^{2}}\geq0\text{
for a.a. }t>0\text{, }x\in\left(  0,1\right)  . \label{CondMax}%
\end{equation}
Hence, $\dfrac{\partial u}{\partial t}-\dfrac{\partial^{2}u}{\partial x^{2}%
}\geq0$ in the sense of distributions as well.
\end{lemma}

\begin{proof}
The function $u$ is the unique solution to problem (\ref{Aux}) for some $f\in
L_{loc}^{2}(\tau,+\infty;H)$ such that $f(t,x)\in b(t)H_{0}(u(t,x))+\omega
(t)u(t,x)$ for a.a. $(t,x)$. If we prove for a.a. $t>0$ that $f\left(
t,x\right)  \geq0$ for a.a. $x\in\left(  0,1\right)  $, then (\ref{CondMax}) follows.

We fix $t>0$ such that the derivative $\dfrac{du}{dt}\left(  t\right)  $
exists in the classical sense and $f\left(  t,\text{\textperiodcentered
}\right)  \in H$. Denote $A_{t}=\{x\in\lbrack0,1]:u\left(  t,x\right)  >0\}$,
$B_{t}=\{x\in\lbrack0,1]:u\left(  t,x\right)  =0\},\ A_{qt}=A_{t}%
\cap\mathbb{Q}$. Since $u\left(  t\right)  \in V\subset C([0,1])$, for any
$x\in A_{qt}$ there exists a maximal interval $I_{x}=[x_{1},x_{2}%
],\ x_{1}<x_{2}$, $x\in I_{x}$, such that $u\left(  y\right)  >0,$ for any
$y\in\left(  x_{1},x_{2}\right)  ,$ and $u\left(  x_{1}\right)  =u\left(
x_{2}\right)  =0$. Again by continuity, for any $x\in A_{t}$ there exists
$\overline{x}\in A_{qt}$ such that $x\in int\ I_{\overline{x}}$. Hence,
$A_{t}=\cup_{x\in A_{qt}}int\ I_{x}$, so $A_{t}$ is the countable union of
open intervals $I_{i}=(x_{i},x_{i+1}),$ where $u\left(  t,x\right)  >0$ for
$x\in\left(  x_{i},x_{i+1}\right)  $ and $u\left(  t,x_{i}\right)  =u\left(
t,x_{i+1}\right)  =0$. Therefore, $B_{t}$ is the countable union of closed
intervals $\widetilde{I}_{j}=[\widetilde{x}_{j},\widetilde{x}_{j+1}],$ where
$u\left(  t,x\right)  =0$ for $x\in\widetilde{I}_{j}$. It is possible in this
last case that $x_{i}=x_{i+1}$. We split $B_{t}$ into two sets:%
\[
B_{t}=B_{t}^{1}\cup B_{t}^{2},
\]%
\begin{align*}
B_{t}^{1}  &  =\{x\in B_{t}:x\in(\widetilde{x}_{i},\widetilde{x}%
_{i+1}),\ \widetilde{x}_{i}<\widetilde{x}_{i+1},\ u\left(  t,y\right)
=0\ \forall y\in\lbrack\widetilde{x}_{i},\widetilde{x}_{i+1}]\},\\
B_{t}^{2}  &  =\{x\in B_{t}:u\left(  t,x\right)  =0,\ \text{there is no
}\varepsilon>0\text{ such that }u\left(  t,y\right)  =0\text{ for }y\in\lbrack
x-\varepsilon,x+\varepsilon]\}.
\end{align*}
The set $B_{t}^{2}$ has zero measure. Thus, it is enough to check that
$f\left(  t,x\right)  \geq0$ for a.a. $x\in A_{t}\cup B_{t}^{1}.$

If $x\in A_{t}$, then $H_{0}(u(t,x))=1$, so that $f\left(  t,x\right)
=b\left(  t\right)  +\omega\left(  t\right)  u\left(  t,x\right)  >0.$

If $x\in(\widetilde{x}_{i},\widetilde{x}_{i+1})$, where $(\widetilde{x}%
_{i},\widetilde{x}_{i+1})$ is an interval given in the definition of
$B_{t}^{1}$, then the derivative $\dfrac{\partial^{2}u}{\partial x^{2}}\left(
t,x\right)  $ exists in the classical sense and is equal to zero. This means
that
\[
\frac{\partial u}{\partial t}\left(  t,x\right)  =f\left(  t,x\right)  \text{
for a.a. }x\in(\widetilde{x}_{i},\widetilde{x}_{i+1}).
\]
Since $\dfrac{du}{dt}\left(  t\right)  $ exists in the classical sense,
$\dfrac{u\left(  t+h\right)  -u\left(  t\right)  }{h}-\dfrac{du}{dt}\left(
t\right)  \rightarrow0$ in $H$, so $\dfrac{u\left(  t+h_{n},x\right)
-u\left(  t,x\right)  }{h_{n}}-\dfrac{\partial u}{\partial t}\left(
t,x\right)  \rightarrow0$ for a.a. $x\in\left(  0,1\right)  $ and some
subsequence. For a.a. $x\in(\widetilde{x}_{i},\widetilde{x}_{i+1})$ we have
that $\dfrac{u\left(  t+h,x\right)  -u\left(  t,x\right)  }{h}\geq0$ (as
$u\left(  t+h,x\right)  \geq0$), and then $\dfrac{\partial u}{\partial
t}\left(  t,x\right)  \geq0$ for a.a. $x\in(\widetilde{x}_{i},\widetilde{x}%
_{i+1})$. Thus, $f\left(  t,x\right)  \geq0$ for a.a. $x\in(\widetilde{x}%
_{i},\widetilde{x}_{i+1}).$
\end{proof}

\bigskip

\begin{lemma}
\label{PositiveSolutions}Assume that (\ref{Cond1}) holds. Let $u_{\tau}\in H$
be such that $u_{\tau}\geq0$ but $u_{\tau}\not \equiv 0$ and let $u\left(
\text{\textperiodcentered}\right)  $ be a non-negative solution to problem
(\ref{Incl}). Then the solution $u\left(  \text{\textperiodcentered}\right)  $
is unique in the class of non-negative solutions and $u(t)$ is positive for
any $t>\tau.$
\end{lemma}

\begin{proof}
We will use the minimum principle for non-smooth functions proved in
\cite{Kadlec} (see the appendix) in order to prove that $u\left(  t\right)  $
is positive for any $t>\tau.$ By contradiction we assume the existence of
$t_{0}>\tau$ and $x_{0}\in\left(  0,1\right)  $ such that $u\left(
t_{0},x_{0}\right)  =0$. Let
\begin{align*}
\mathcal{O}  &  \mathcal{=[}\tau,t_{0}]\times\lbrack0,1],\\
Q_{\rho,\sigma}  &  =\{(t,x):t\in(t_{0}-\sigma,t_{0}),\left\vert x-\frac{1}%
{2}\right\vert <\rho\},
\end{align*}
where $\sigma=t_{0}-\tau$. Choosing $\rho=\frac{1}{2}$ we have that%
\begin{equation}
\underset{(t,x)\in Q_{\nu\rho,\sigma_{1}}}{\inf}\ u(t,x)=0, \label{EssInf}%
\end{equation}
for some $0<\nu<1$ and any $0<\sigma_{1}<\sigma$. By (\ref{CondMax}) and
Theorem \ref{Minimum} we obtain that $u(t,x)=0$ for a.a. $\left(  t,x\right)
\in Q_{\frac{1}{2},\sigma}$, which is not possible.

It remains to prove the uniqueness. Since $u\left(  t\right)  $ is positive
for $t>\tau$, $u$ has to be the solution of problem (\ref{Aux}) with
$f(t,x)=b(t)+\omega(t)u(t,x)$. Let there exist two solutions $u_{1},u_{2}$.
Then the difference $v=u_{1}-u_{2}$ satisfies the equality%
\[
\frac{\partial v}{\partial t}-\dfrac{\partial^{2}v}{\partial x^{2}}%
=\omega(t)v.
\]
Multiplying by $v$ and making use of Gronwall's lemma we get%
\[
\left\Vert v(t)\right\Vert ^{2}\leq e^{2\int_{\tau}^{t}\omega\left(  s\right)
ds}\left\Vert v(\tau)\right\Vert ^{2}=0.
\]

\end{proof}

\bigskip

\begin{lemma}
\label{ZeroSolutions}Assume that (\ref{Cond1}) holds. If $u_{\tau}\equiv0$,
then, apart from the zero solution, all the possible non-negative solutions
are of the following type:%
\begin{equation}
u(t)=\left\{
\begin{array}
[c]{c}%
0\text{ if }\tau\leq t\leq t_{0},\\
u_{t_{0}}(t)\text{ if }t\geq t_{0},
\end{array}
\right.  \label{Sol0}%
\end{equation}
where $u_{t_{0}}($\textperiodcentered$)$ is the unique solution to the problem%
\begin{equation}
\left\{
\begin{array}
[c]{l}%
\dfrac{\partial u}{\partial t}-\dfrac{\partial^{2}u}{\partial x^{2}%
}=b(t)+\omega(t)u(t),\text{ on }(t_{0},\infty)\times(0,1),\\
u(t,0)=u(t,1)=0,\\
u(t_{0},x)=0,
\end{array}
\right.  \label{Aux3}%
\end{equation}
and $u(t)$ is positive for all $t>t_{0}.$
\end{lemma}

\begin{proof}
By comparison, it is clear that the function given in (\ref{Sol0}) is a
non-negative solution to problem (\ref{Incl}) for any $t_{0}\geq\tau$. It
remains to check that these are the only ones and that $u(t)$ is positive for
$t>t_{0}$.

Since any solution $u\left(  \text{\textperiodcentered}\right)  $ to problem
(\ref{Incl}) with $u_{\tau}\equiv0$ satisfies that $u\in C([\tau
,+\infty),C([0,1])$, if $u\left(  \overline{t}\right)  $ is not identically
equal to $0$, then there exists an interval $\left(  t_{0},t_{1}\right)  $
containing $\overline{t}$ such that $u(t)$ is not identically equal to $0$ for
all $t\in\left(  t_{0},t_{1}\right)  $. By Lemma \ref{PositiveSolutions} the
element $u(t)$ is positive for any $t>t_{0}$. Thus, if $u\left(  \overline
{t}\right)  \not \equiv 0$ for some $\overline{t}>\tau$, then $u\left(
t\right)  $ is positive for any $t\geq\overline{t}$. This implies that any
non-negative solution not being identically equal to $0$ has to be of the type
(\ref{Sol0}) and $u\left(  t\right)  $ is positive for any $t>t_{0}$.
Obviously, $u_{t_{0}}($\textperiodcentered$)$ has to be the solution to
problem (\ref{Aux3}).
\end{proof}

\bigskip

\begin{corollary}
\label{Backward}Assume that (\ref{Cond1}) holds. If $u_{\tau}\in V$ be such
that $u_{\tau}\geq0$, $u_{\tau}\left(  x_{0}\right)  =0$ at some $x_{0}%
\in\left(  0,1\right)  $ but $u_{\tau}\not \equiv 0$, then there cannot exist
a non-negative solution backwards in time.
\end{corollary}

\begin{proof}
By contradiction let there be a non-negative solution $u\left(
\text{\textperiodcentered}\right)  $ to problem (\ref{Incl}) in some interval
$[t_{0},+\infty)$ with $t_{0}<\tau$ such that $u\left(  \tau\right)  =u_{\tau
}$. By Lemmas \ref{PositiveSolutions} and \ref{ZeroSolutions} we have that
$u\left(  \tau\right)  $ needs to be positive, so $u_{\tau}\left(
x_{0}\right)  >0$.
\end{proof}

\begin{remark}
Assume that (\ref{Cond1}) holds. If $u_{\tau}\in H$ but $u_{\tau}\not \in V$,
then there cannot exist a non-negative solution $u\left(
\text{\textperiodcentered}\right)  $ backwards in time as well. This follows
from the fact that $u\in C(\left(  \tau,+\infty\right)  ,V).$
\end{remark}

\begin{corollary}
\label{Backward2}Assume that (\ref{Cond1}) holds. If $u_{\tau}\equiv0$, then
the unique non-negative solution backwards in time is the zero solution, that
is, $u\left(  t\right)  \equiv0$ for $t\leq\tau$.
\end{corollary}

\bigskip

Adding some extra conditions on the functions $b,\omega,$ we will establish
that the solutions are more regular. First, we state a technical lemma about
the derivatives of the product of two functions.

\begin{lemma}
\label{LemaProducto}Let $S=\left(  T_{1},T_{2}\right)  \subset\mathbb{R},$
$h\in W^{1,2}(S)$, $v\in W^{1,2}(S;H)$. Then $hv\in W^{1,2}(S;H)$ and
\begin{equation}
\frac{d}{dt}\left(  hv\right)  =h\frac{dv}{dt}+\frac{dh}{dt}v.
\label{DerivProduct}%
\end{equation}

If $h\in W^{k,2}(S)$, $v\in W^{k,2}(S;H),$ $k\in\mathbb{N}$, then $hv\in
W^{k,2}(S;H)$ and%
\begin{equation}
\frac{d^{j}}{dt^{j}}\left(  hv\right)  =\sum_{i=0}^{j}a_{i}\frac{d^{i}%
h}{dt^{i}}\frac{d^{j-i}v}{dt^{j-i}},\text{ for }j=1,...,k,
\label{DerivProduct2}%
\end{equation}
where $a_{i}$ are the coefficients of the Tartaglia/Pascal triangle.
\end{lemma}

\begin{proof}
Since $W^{1,2}(S)\subset L^{\infty}(S)$, we get that $hv\in L^{2}(S;H)$. Let
$h_{n}\in C^{1}(\overline{S})$, $v_{n}\in C^{1}(\overline{S},H)$ be such that
$h_{n}\rightarrow h$ in $W^{1,2}(S)$, $v_{n}\rightarrow v$ in $W^{1,2}(S;H)$
(see \cite[Chapter IV]{Gajewsky}). Then%
\[
\frac{d}{dt}\left(  h_{n}v_{n}\right)  =h_{n}\frac{dv_{n}}{dt}+\frac{dh_{n}%
}{dt}v_{n}\rightarrow h\frac{dv}{dt}+\frac{dh}{dt}v\text{ in }L^{1}(S;H).
\]
Thus, (\ref{DerivProduct}) follows and $h\in L^{\infty}(S)$, $v\in
C(\overline{S},H)\subset L^{\infty}(S;H)$ imply that $hv\in W^{1,2}(S;H).$

Let $h\in W^{k,2}(S)$, $v\in W^{k,2}(S;H),$ $k\in\mathbb{N}$. By induction
assume that $h\in W^{m,2}(S)$, $v\in W^{m,2}(S;H)$ and (\ref{DerivProduct2})
holds for $m\in\{1,...,j\}$ with $j\in\{1,...,k-1\}$. Since $\dfrac{d^{i}%
h}{dt^{i}}\in W^{1,2}(S)$, $\dfrac{d^{j-i}v}{dt^{j-i}}\in W^{1,2}(S;H)$ for
any $i\in\{0,...,j\}$, the previous result gives that $\dfrac{d^{i}h}{dt^{i}%
}\dfrac{d^{j-i}v}{dt^{j-i}}\in W^{1,2}(S;H)$ and
\begin{equation}
\frac{d}{dt}\left(  \dfrac{d^{i}h}{dt^{i}}\dfrac{d^{j-i}v}{dt^{j-i}}\right)
=\dfrac{d^{i+1}h}{dt^{i+1}}\dfrac{d^{j-i}v}{dt^{j-i}}+\dfrac{d^{i}h}{dt^{i}%
}\dfrac{d^{j-i+1}v}{dt^{j-i+1}}. \label{Sum1}%
\end{equation}
Hence, $hv\in W^{j+1,2}(S;H)$. For any $i\in\{1,...,j\}$ we have%
\begin{equation}
\frac{d}{dt}\left(  \dfrac{d^{i-1}h}{dt^{i-1}}\dfrac{d^{j-i+1}v}{dt^{j-i+1}%
}\right)  =\dfrac{d^{i}h}{dt^{i}}\dfrac{d^{j-i+1}v}{dt^{j-i+1}}+\dfrac
{d^{i-1}h}{dt^{i-1}}\dfrac{d^{j-i+2}v}{dt^{j-i+2}}, \label{Sum2}%
\end{equation}
so the second term in (\ref{Sum1}) and the first term in (\ref{Sum2}) are
equal. Thus,%
\[
\frac{d^{j+1}}{dt^{j+1}}\left(  hv\right)  =a_{0}h\frac{d^{j+1}v}{dt^{j+1}%
}+\sum_{i=1}^{j}\left(  a_{i-1}+a_{i}\right)  \dfrac{d^{i}h}{dt^{i}}%
\dfrac{d^{j-i+1}v}{dt^{j-i+1}}+a_{j}\dfrac{d^{j+1}h}{dt^{j+1}}v,
\]
proving formula (\ref{DerivProduct2}) for $j+1$.
\end{proof}

\begin{lemma}
\label{H3Regularity}Assume that (\ref{Cond1}) holds and, additionally, that
$b,\omega\in W_{loc}^{1,2}(\mathbb{R})$. Let $u_{\tau}\in H$ be such that
$u_{\tau}\geq0$ and $\left\Vert u_{\tau}\right\Vert >0$. Then the unique
non-negative solution $u\left(  \text{\textperiodcentered}\right)  $ to
problem (\ref{Incl}) satisfies that
\begin{align}
u  &  \in C((\tau,+\infty),H^{3}(0,1))\cap C^{1}((\tau,+\infty
),V),\label{URegularity}\\
\dfrac{d^{2}u}{dt^{2}}  &  \in L_{loc}^{2}(\tau+\varepsilon,+\infty;H)\text{
for all }\varepsilon>0.\nonumber
\end{align}
The partial derivatives $u_{t},u_{xx}$ exist in the classical sense and are
continuous on $(\tau,+\infty)\times\lbrack0,1]$.

If $b,\omega\in W_{loc}^{k+1,2}(\mathbb{R})$, where $k\in\mathbb{N}$, then%
\begin{align}
u  &  \in C^{k+1}(\left(  \tau,+\infty\right)  ,V),\ u\in\cap_{j=0}^{k}%
C^{j}\left(  (\tau,+\infty),H^{2(k-j)+3}(0,1)\right)  ,\label{URegularity2}\\
\dfrac{d^{k+2}u}{dt^{k+2}}  &  \in L_{loc}^{2}(\tau+\varepsilon,+\infty
;H)\text{ for all }\varepsilon>0,\ u\in C^{k+1}(\left(  \tau,+\infty\right)
\times\lbrack0,1]).\nonumber
\end{align}

\end{lemma}

\begin{proof}
By definition $u$ is the unique solution of problem (\ref{Aux})\ with
$f(t,x)\in b(t)H_{0}(u(t,x))+\omega(t)u(t,x)$ for a.a. $(t,x)$. In view of
Lemma \ref{PositiveSolutions} $u(t)$ is positive for any $t>\tau$, so
$H_{0}(u(t,x))=1$ for a.a. $(t,x)\in\left(  \tau,+\infty\right)  \times\left(
0,1\right)  $. Since $u\in C((\tau,+\infty),V)$, $\dfrac{du}{dt}\in
L_{loc}^{2}(\tau+\varepsilon,+\infty;H)$, for all $\varepsilon>0$, and
$b,\omega\in W_{loc}^{1,2}(\mathbb{R})$, we obtain by Lemma \ref{LemaProducto}
that
\begin{equation}
f\in L_{loc}^{2}(\tau,+\infty;H)\cap C((\tau,+\infty),H^{1}(0,1))\cap
W_{loc}^{1,2}(\tau+\varepsilon,+\infty;H)\text{ for all }\varepsilon>0.
\label{fRegularity}%
\end{equation}
The first statement follows from Proposition \ref{Regularity} and Corollary
\ref{Regularity2} (see the appendix).

Let now $b,\omega\in W_{loc}^{k+1,2}(\mathbb{R})$ with $k\in\mathbb{N}$. By
induction assume that for $j\in\{1,...,k\}$ we have that $u\in C^{j}(\left(
\tau,+\infty\right)  ,V)$ and
\begin{align*}
u  &  \in\cap_{i=0}^{j-1}C^{i}\left(  (\tau,+\infty),H^{2(j-i)+1}(0,1)\right)
,\\
\dfrac{d^{j+1}u}{dt^{j+1}}  &  \in L_{loc}^{2}(\tau+\varepsilon,+\infty
;H)\text{ for all }\varepsilon>0.
\end{align*}
By Lemma \ref{LemaProducto} we get $f\in W_{loc}^{j+1,2}(\tau+\varepsilon
,+\infty;H),$ for all $\varepsilon>0$, and%
\[
f\in\cap_{i=0}^{j}C^{i}\left(  (\tau,+\infty),H^{2(j-i)+1}(0,1)\right)  .
\]
Now, Lemma \ref{Regularity3} and Corollary \ref{Regularity4} yield%
\begin{align*}
\dfrac{d^{j+2}u}{dt^{j+2}}  &  \in L_{loc}^{2}(\tau+\varepsilon,+\infty
;H)\text{ for all }\varepsilon>0,\ u\in C^{j+1}(\left(  \tau,+\infty\right)
,V),\\
u  &  \in\cap_{i=0}^{j}C^{i}\left(  (\tau,+\infty),H^{2(j-i)+3}(0,1)\right)
,\ u\in C^{j+1}(\left(  \tau,+\infty\right)  \times\lbrack0,1]).
\end{align*}

\end{proof}

\bigskip

Adding some extra conditions on the initial condition and the functions
$b,\omega,$ the solution will now be proved to be unique in the class of all solutions.

Let $V^{2r}=D(A^{r})$. Take $u_{\tau}\in V^{2r}$ with $\frac{3}{4}<r<1$.
Hence, $u\in C([\tau,+\infty),V)$ implies that $u$ is a solution of problem
(\ref{Aux}) with $f\in L_{loc}^{\infty}(\tau,+\infty;L^{\infty}(0,1))$, so
Lemma 42.7 in \cite{SellBook} implies that
\[
u\in C([\tau,+\infty),V^{2r}).
\]
We observe that $H^{s}\left(  0,1\right)  $ is continuously embedded in
$C^{1}([0,1])$ if $s>\frac{3}{2}$ \cite[Lemma 4.4]{Grisvard}. This fact,
together with $D(A^{r})=H^{2r}(0,1)\cap H_{0}^{1}\left(  0,1\right)  $ for
$\frac{3}{4}<r<1$ \cite[Theorem 1]{Fujiwara}, implies that
\begin{equation}
u\in C([\tau,+\infty),C^{1}([0,1])). \label{ContDeriv}%
\end{equation}

\begin{lemma}
\label{PositiveSolutions3}Let $u_{\tau}\in V^{2r},\ \frac{3}{4}<r<1,$ be
positive and such that $\frac{d}{dx}u_{\tau}(0)>0,\ \frac{d}{dx}u_{\tau}%
(1)<0$. Assume that (\ref{Cond1}) holds and, additionally, that $b,\omega\in
W_{loc}^{1,2}(\mathbb{R})$. Then there is a unique solution $u\left(
\text{\textperiodcentered}\right)  $ to problem (\ref{Incl}), which satisfies
that $u(t)$ is positive for any $t\geq\tau$ and (\ref{URegularity}) as well.
Moreover, $u_{x}(t,0)>0,\ u_{x}(t,1)<0$ for any $t\geq\tau$ and the partial
derivatives $u_{t},u_{xx}$ exists in the classical sense and are continuous on
$(\tau,+\infty)\times\lbrack0,1]$.

If, moreover, $b,\omega\in W_{loc}^{k+1,2}(\mathbb{R})$, , where
$k\in\mathbb{N}$, then additionally $u$ satisfies (\ref{URegularity2}).
\end{lemma}

\begin{proof}
Let $u\left(  \text{\textperiodcentered}\right)  $ be an arbitrary solution to
problem (\ref{Incl}). First we will show that there exists an interval
$[\tau,t_{1}]$, $t_{1}>\tau$, such that $u(t)$ is positive and $u_{x}%
(t,0)>0,\ u_{x}(t,1)<0,\ $for any $t\in\lbrack\tau,t_{1}]$. Since $u_{x}$ is
jointly continuous by (\ref{ContDeriv}), there exist $t_{1}>\tau
,\ 0<x_{0}<x_{1}<1$ and $\alpha_{0}>0$ such that
\begin{align*}
u_{x}(t,x)  &  \geq\alpha_{0},\ \forall\ t\in\lbrack\tau,t_{1}],\ x\in
\lbrack0,x_{0}]\\
u_{x}(t,1)  &  \leq-\alpha_{0}\text{, }\forall\ t\in\lbrack\tau,t_{1}%
],\ x\in\lbrack x_{1},1].
\end{align*}
Hence,
\begin{align*}
u\left(  t,x\right)   &  \geq\alpha_{0}x\text{ for }x\in\lbrack0,x_{0}%
],\ t\in\lbrack\tau,t_{1}],\\
u(t,x)  &  \geq\alpha_{0}(1-x)\ \text{for }x\in\lbrack x_{1},1],\ t\in
\lbrack\tau,t_{1}].
\end{align*}
Finally, by the joint continuity of $u$ there exist $\tau<t_{2}\leq t_{1}$ and
$\alpha_{1}>0$ such that%
\[
u(t,x)\geq\alpha_{1}\text{ for }x\in\lbrack x_{0},x_{1}],\ t\in\lbrack
\tau,t_{2}].
\]

We state that in fact these properties hold in any interval $[\tau,t]$.

We observe first that there cannot be a time $t_{0}>\tau$ and $x_{0}\in\left(
0,1\right)  $ such that%
\begin{align}
u\left(  t,x\right)   &  >0\text{ }\forall\ \tau\leq t<t_{0},\ x\in
\lbrack0,1],\label{uPositive}\\
u(t_{0},x)  &  \geq0,\ \forall x\in\lbrack0,1],\ u\left(  t_{0},x_{0}\right)
=0.\nonumber
\end{align}
In such a case, $u(t_{0},x_{0})$ is the minimum of the function $u$ in the
region $\Omega=\{(t,x):\tau<t\leq t_{0},\ 0\leq x\leq1\}$.

Since by Lemma \ref{PositiveSolutions} $u\left(  \text{\textperiodcentered
}\right)  $ is the unique non-negative solution in the interval $[\tau,t_{0}%
]$, Lemma \ref{H3Regularity} implies that the derivatives $u_{t},u_{xx}$ exist
in the classical sense and are continuous on $\Omega$. Then, as $u(t_{0}%
,x_{0})$ is the minimum in $\Omega$, we have%
\[
u_{t}(t_{0},x_{0})\leq0,\ u_{xx}(t_{0},x_{0})\geq0.
\]
Since%
\[
u_{t}(t,x)-u_{xx}(t,x)=b(t)+\omega(t)u(t,x)\text{ for }t\in\left(  \tau
,t_{0}\right)  ,\ x\in\left(  0,1\right)  ,
\]
we obtain by continuity that%
\[
0\geq u_{t}(t_{0},x_{0})-u_{xx}(t_{0},x_{0})=b(t_{0})+\omega(t_{0}%
)u(t_{0},x_{0})>0,
\]
which is a contradiction. Therefore, (\ref{uPositive})\ cannot happen.

Further, we will prove that while $u\left(  t\right)  $ is positive the
spatial derivatives at the boundary cannot vanish. By contradiction assume for
example that for some $t_{0}>\tau$ one has
\begin{align*}
u_{x}(t_{0},0)  &  =0,\\
u_{x}(t,0)  &  >0,\ \text{for }t\in\lbrack\tau,t_{0}),\\
u(t,x)  &  >0\text{ for }t\in\lbrack\tau,t_{0}],\ x\in\left(  0,1\right)  .
\end{align*}
As we have seen before, the derivatives $u_{t},u_{xx}$ exist in the classical
sense and are continuous on $(\tau,t_{0}]\times\lbrack0,1]$. It is not
possible that $u_{xx}(t_{0},x)<0$ for $x$ in some interval $[0,\varepsilon]$
as in such a case we would have that%
\[
u_{x}(t_{0},x)-u_{x}(t_{0},0)=u_{x}(t_{0},x)=\int_{0}^{x}u_{yy}(t_{0}%
,y)dy<0,\text{ for }x\in(0,\varepsilon],
\]
and then%
\[
u(t_{0},\varepsilon)=\int_{0}^{\varepsilon}u_{y}(t_{0},y)dy<0,
\]
which is false. Therefore, there has to exist a sequence $u_{xx}(t_{0}%
,x_{n})\geq0$ with $x_{n}\rightarrow0^{+}$. Thus,%
\[
u_{t}(t_{0},x_{n})=b(t_{0})+\omega\left(  t_{0}\right)  u\left(  t_{0}%
,x_{n}\right)  +u_{xx}(t_{0},x_{n})\geq b_{0}>0.
\]
By the uniform continuity of $u_{t}$ on the compact set $[\widetilde{\tau
},t_{0}]\times\lbrack0,1],$ where $\widetilde{\tau}>\tau$, there exists
$\delta>0$ such that
\[
u_{t}(t,x_{n})\geq\frac{b_{0}}{2}\text{ for all }n\text{, }t\in\lbrack
t_{0}-\delta,t_{0}].
\]
Hence,%
\[
u(t_{0},x_{n})=u(t_{0}-\delta,x_{n})+\int_{t_{0}-\delta}^{t_{0}}u_{t}%
(t,x_{n})dt\geq\frac{\delta b_{0}}{2}\text{ for all }n,
\]
so%
\[
u(t_{0},x_{n})\rightarrow u(t_{0},0)=0
\]
gives a contradiction. By the same argument $u_{x}$ cannot vanish at $x=1$.

The only possibility left is that $u$ becomes $0$ near the boundary
instantaneously at some moment of time, that is, there exist for example a
time $t_{0}>\tau$ and sequences $t_{n}>t_{0},\ x_{n}>0$ such that
$t_{n}\rightarrow t_{0},\ x_{n}\rightarrow0$ and
\begin{align}
u(t,x)  &  >0\text{ for }t\in\lbrack\tau,t_{0}],\ x\in\left(  0,1\right)
,\label{ZeroBoundary}\\
u(t_{n},x_{n})  &  =0.\nonumber
\end{align}
We have seen that $u_{x}(t,0)>0$ for any $t\in\lbrack\tau,t_{0}]$. By the
joint continuity of $u_{x}$ there exist $t_{1}>t_{0},\ x_{0}>0$ and
$\alpha_{0}>0$ such that $u_{x}\left(  t,x\right)  \geq\alpha_{0}$ for any
$t\in\lbrack t_{0},t_{1}],\ x\in\lbrack0,x_{0}]$. Thus, $u\left(  t,x\right)
\geq\alpha_{0}x$ for $t\in\lbrack t_{0},t_{1}],$ $x\in\lbrack0,x_{0}]$, which
contradicts (\ref{ZeroBoundary}). In the same way one can show that $u$ cannot
become $0$ instantaneously near $x=1$.

We have proved that $u\left(  t\right)  $ is positive and $u_{x}%
(t,0)>0,\ u_{x}(t,1)<0$ for any $t\geq\tau$. Thus, Lemma \ref{H3Regularity}
implies that the solution $u\left(  \text{\textperiodcentered}\right)  $ is
unique in the class of all solutions and the regularity of $u$ as well.
\end{proof}

\bigskip

\section{Structure of the pullback attractor in the positive cone}

In this section, we will apply the previous results in order to study the
structure and regularity of the pullback attractor for the multivalued process
generated by the solutions of problem (\ref{Incl}) in the positive cone.

\subsection{The autonomous case}

We start with the autonomous case, that is, we assume that $b(t)\equiv
b>0,\ \omega(t)\equiv\omega\geq0$. Hence, problem (\ref{Incl}) becomes
\begin{equation}
\left\{
\begin{array}
[c]{l}%
\dfrac{\partial u}{\partial t}-\dfrac{\partial^{2}u}{\partial x^{2}}\in
bH_{0}(u)+\omega u,\text{ on }(0,\infty)\times(0,1),\\
u(t,0)=u(t,1)=0,\\
u(0,x)=u_{0}(x),\text{ }x\in\left(  0,1\right)  .
\end{array}
\right.  \label{InclAut}%
\end{equation}

Additionally, in order to guarantee the existence of the global attractor, we
will suppose throughout this section that%
\begin{equation}
\omega<\pi^{2}, \label{CondAttr}%
\end{equation}
where $\pi^{2}$ is the first eigenvalue of the operator $-\frac{\partial^{2}%
u}{\partial x^{2}}$ in $V$.

Let us denote by $\mathcal{D}(u_{0})$ the set of all solutions of the
autonomous problem (\ref{InclAut}) with initial condition $u_{0}$ at $t=0$.
Let $\mathcal{R}_{0}=\cup_{u_{0}\in H}\mathcal{D}(u_{0})$ be the set of all
solutions. Denoting by $P(H)$ the set of all non-empty subsets of $H$, the
multivalued map $G:\mathbb{R}^{+}\times H\rightarrow P(H)$ is defined by%
\[
G(t,u_{0})=\{u(t):u\in\mathcal{D}(u_{0})\}.
\]
We recall \cite{Valero2001, Valero2005} that this map satisfies that
$G(0,$\textperiodcentered$)$ is the identity map and $G(t+s,u_{0}%
)=G(t,G(s,u_{0}))$ for all $t,s\in\mathbb{R}^{+}$, $u_{0}\in H$, that is, it
is a strict multivalued semiflow, and that it possesses a global compact
connected invariant attractor $\mathcal{A}$ in the phase space $H$, which
means that:

\begin{itemize}
\item $\mathcal{A}$ is compact in $H;$

\item $\mathcal{A}$ is connected in $H;$

\item $\mathcal{A=}G(t,\mathcal{A})$ for all $t\in\mathbb{R}^{+}$ (that is, it
is strictly invariant);

\item $dist_{H}(G(t,B),\mathcal{A})\rightarrow0$ as $t\rightarrow+\infty$ for
any set $B$ bounded in $H$ (that is, it is attracting).
\end{itemize}

The function $\varphi:\mathbb{R}\rightarrow H$ is said to be a complete
trajectory of $\mathcal{R}_{0}$ if $\varphi\left(  \text{\textperiodcentered
}+h\right)  \mid_{t\geq0}\in\mathcal{R}_{0}$ for any $h\in\mathbb{R}$. It is
bounded if the set $\cup_{t\in\mathbb{R}}\varphi\left(  t\right)  $ is
bounded. The global attractor $\mathcal{A}$ consists of the union of the
elements of all the bounded complete trajectories \cite{ARV06}, that is,%
\[
\mathcal{A}=\{\varphi(0):\varphi\text{ is a bounded complete trajectory of
}\mathcal{R}_{0}\}.
\]
Moreover, each bounded complete trajectory satisfies that%
\begin{align*}
\varphi(t)  &  \rightarrow z_{1}\text{ as }t\rightarrow+\infty,\\
\varphi(t)  &  \rightarrow z_{2}\text{ as }t\rightarrow-\infty,
\end{align*}
where $z_{i}$ are fixed points such that $E(z_{1})<E(z_{2})$ for the Lyapunov
energy functional $E:V\rightarrow\mathbb{R}$ given by $E(u)=\frac{1}{2}%
\int_{0}^{1}\left\vert \dfrac{du}{dx}\right\vert ^{2}dx-\int_{0}^{1}\left\vert
u\right\vert dx$. Also, the stationary points were fully described in
\cite{ARV06}, where it was shown that there is an infinite but countable
number of them, denoted by $v_{0}=0,\ v_{1}^{\pm},\ v_{2}^{\pm},...,v_{n}%
^{\pm},...$, and that they are ordered by the functional $E$, that is,
$E\left(  0\right)  <E\left(  v_{1}^{\pm}\right)  <E\left(  v_{2}^{\pm
}\right)  <...$ The functions $v_{n}^{\pm}$ have exactly $n-1$ zeros in
$\left(  0,1\right)  $ and $v_{n}^{+}=-v_{n}^{-}$. In particular, $v_{1}^{\pm
}$ have no zeros in $\left(  0,1\right)  $ and they are given by%
\begin{align*}
v_{1}^{+}(x)  &  =\frac{b}{\omega}\cos\left(  \sqrt{\omega}x\right)
+\frac{b\left(  1-\cos\sqrt{\omega}\right)  }{\omega\sin\left(  \sqrt{\omega
}\right)  }\sin\left(  \sqrt{\omega}x\right)  -\frac{b}{\omega}\text{ if
}0<\omega<\pi^{2},\\
v_{1}^{+}(x)  &  =-\frac{x^{2}}{2}+\frac{x}{2}\text{ if }\omega=0.
\end{align*}
It is clear that $v_{1}^{\pm}\in C^{\infty}\left(  [0,1]\right)  $. By Lemma
\ref{PositiveSolutions3} the solutions starting at $v_{1}^{\pm}$ are unique.
In addition, we know from \cite[Section 3]{CVL19} that
\[
v_{1}^{-}\leq z\leq v_{1}^{+}\text{ for any }z\in\mathcal{A}.
\]

In relation with the regularity of the attractor, it was proved in
\cite[Theorem 3.1 and Corollary 3.2]{ARV06} that the global attractor is
compact in $W^{2-\delta,p}(0,1)$ for any $\delta>0$, $1\leq p<\infty$ and that
it is bounded in $L^{\infty}(0,1)$. In our particular one dimensional
situation, the last result is also consequence of the continuous embedding
$H_{0}^{1}(0,1)\subset L^{\infty}(0,1).$ Let us extend these results to the
$V^{2r}$ spaces.

\begin{lemma}
\label{CompactV2r}Let (\ref{CondAttr}) hold. The global attractor
$\mathcal{A}$ is compact in $V^{2r}$ for any $0\leq r<1$.
\end{lemma}

\begin{proof}
For any $z\in\mathcal{A}$ the invariance of $\mathcal{A}$ implies the
existence of $u\in\mathcal{D}_{0}(u_{0})$, $u_{0}\in\mathcal{A}$, such that
$z=u(1)$ and $u\left(  t\right)  \in\mathcal{A}$ for all $t\geq0.$ As $u$ is a
mild solution, by the variation of constants formula we get%
\[
z=e^{-A}u_{0}+\int_{0}^{1}e^{-A(1-s)}f(s)ds,
\]
where $f(s,x)\in bH_{0}(u(t,x))+\omega u(t,x)$ for a.a. $(t,x)$. By the
boundedness of $\mathcal{A}$ in $H,$ there exists a universal constant such
that
\[
\left\Vert u_{0}\right\Vert \leq C,\ \left\Vert f\right\Vert _{L^{\infty
}(0,1;H)}\leq C,
\]
where $C$ does not depend on the chosen $z$. Then by standard estimates of the
norm of $e^{-At}$ in $V^{2r}$ \cite[Theorem 37.5]{SellBook} for some constants
$M_{r}>0$, $a\in\mathbb{R}$ we have%
\begin{align*}
\left\Vert A^{r}z\right\Vert  &  \leq\left\Vert A^{r}e^{-A}u_{0}\right\Vert
+\int_{0}^{1}\left\Vert A^{r}e^{-A(1-s)}f(s)\right\Vert ds\\
&  \leq M_{r}e^{-a}C+M_{r}C\int_{0}^{1}\left(  1-s\right)  ^{-r}ds\leq C_{r}.
\end{align*}
Therefore, $\mathcal{A}$ is bounded in $V^{2r}$ if $r<1$.

The compact embedding $V^{\alpha}\subset V^{\beta}$ for $\alpha>\beta$ implies
that $\mathcal{A}$ is relatively compact in $V^{2r}$ for any $r<1$. Since
$\mathcal{A}$ is closed in $H,$ it is closed in $V^{2r}$ as well, so
$\mathcal{A}$ is compact in $V^{2r}$ for any $0\leq r<1$.
\end{proof}

\bigskip

Let us consider non-negative solutions to problem (\ref{Incl}). We have seen
that for any $u_{0}\in H$ satisfying $u_{0}\geq0$ there exists at least one
solution $u\left(  \text{\textperiodcentered}\right)  $ such that $u\left(
t\right)  \geq0$ for any $t\geq0$. We denote by $\mathcal{D}^{+}(u_{0})$ the
set of all non-negative solutions of problem (\ref{Incl}) with initial
condition $u_{0}\geq0$ at time $t=0$ and let $\mathcal{R}_{0}^{+}=\cup
_{u_{0}\in H}\mathcal{D}^{+}(u_{0})$. Also, let $H^{+}$ be the positive cone
of $H$, that is,%
\[
H^{+}=\{v\in H:v(x)\geq0\text{ for a.a. }x\in\left(  0,1\right)  \}.
\]

We define the strict multivalued semiflow $G^{+}:\mathbb{R}^{+}\times
H^{+}\rightarrow P(H^{+})$ given by%
\[
G^{+}(t,u_{0})=\{u(t):u\in\mathcal{D}_{0}^{+}(u_{0})\}.
\]
The bounded complete trajectories of $\mathcal{R}_{0}^{+}$ are all the bounded
complete trajectories of $\mathcal{R}_{0}$ such that $\varphi(t)\geq0$ for any
$t\in\mathbb{R}$. The semiflow $G^{+}$ possesses a global compact connected
invariant attractor, denoted by $\mathcal{A}^{+}$ and%
\[
\mathcal{A}^{+}=\{\varphi(0):\varphi\text{ is a bounded complete trajectory of
}\mathcal{R}_{0}^{+}\}.
\]
In the positive cone the only fixed points are $0$ and $v_{1}^{+}$. Therefore,
any bounded complete trajectory different from $0$ or $v_{1}^{+}$ has to
satisfy that%
\begin{align}
\varphi(t) &  \rightarrow v_{1}^{+}\text{ as }t\rightarrow+\infty
,\label{Connection}\\
\varphi(t) &  \rightarrow0\text{ as }t\rightarrow-\infty.\nonumber
\end{align}
Also,
\[
0\leq z\leq v_{1}^{+}\text{ for any }z\in\mathcal{A}^{+}.
\]

We observe that by Lemma \ref{CompactV2r} the convergence in (\ref{Connection}%
) is true in the $V^{2r}$ spaces with $r<1$ and by the compact embedding
$V^{2r}\subset C^{1}([0,1])$, $r>\frac{3}{4},$ in the space $C^{1}([0,1])$ as well.

Let us consider the set $D=\{v\in V:v(x)>0$ for $x\in\left(  0,1\right)  \}$.

\begin{lemma}
\label{UniqueNondegenerate}Let (\ref{CondAttr}) hold. The fixed point
$v_{1}^{+}$ is the unique bounded complete trajectory $\varphi\left(
\text{\textperiodcentered}\right)  $ for which there exists a time $t_{0}$
satisfying that $\varphi(t)\in D$ for any $t\leq t_{0}.$
\end{lemma}

\begin{proof}
Let there exist another bounded complete trajectory $\varphi\left(
\text{\textperiodcentered}\right)  $ for which there exists a time $t_{0}$
satisfying that $\varphi(t)\in D$ for any $t\leq t_{0}.$ Then they are
solutions in any interval $[\tau,t_{0}]$ of the problem%
\begin{equation}
\left\{
\begin{array}
[c]{l}%
\dfrac{\partial u}{\partial t}-\dfrac{\partial^{2}u}{\partial x^{2}}\in
b+\omega u,\text{ on }(\tau,t_{0})\times(0,1),\\
u|_{\partial\Omega}=0,\\
u(\tau,x)=u_{\tau}(x),
\end{array}
\right.  \label{LinearAut}%
\end{equation}
with $u_{\tau}=v_{1}^{+}$ and $u_{\tau}=\varphi\left(  \tau\right)  $,
respectively. Taking the difference of the two equations and multiplying by
$v=\varphi-v_{1}^{+}$ we obtain that%
\[
\frac{1}{2}\frac{d}{dt}\left\Vert v\right\Vert ^{2}+\left(  \pi^{2}%
-\omega\right)  \left\Vert v\right\Vert ^{2}\leq0,
\]
so%
\begin{equation}
\left\Vert v(t)\right\Vert ^{2}\leq e^{-2\left(  \pi^{2}-\omega\right)
(t-\tau)}\left\Vert v(\tau)\right\Vert ^{2}\rightarrow0\text{ as }%
\tau\rightarrow-\infty\text{ for all }t\leq t_{0}. \label{ConvergBackwards}%
\end{equation}
Hence, $\varphi\left(  t\right)  =v_{1}^{+}$ for all $t\leq t_{0}$. Since the
solution of the autonomous problem (\ref{Incl}) with $u_{\tau}=v_{1}^{+}$ is
unique by Lemma \ref{PositiveSolutions3}, we have that $\varphi\equiv
v_{1}^{+}.$
\end{proof}

\begin{corollary}
\label{Backward3}Let (\ref{CondAttr}) hold. If $\varphi$ is a bounded complete
trajectory such that $\varphi\left(  t\right)  \rightarrow0$ as $t\rightarrow
-\infty$, then for any $t_{1}$ there exists a time $t_{0}\leq t_{1}$ such that
$\varphi\left(  t_{0}\right)  \not \in D.$
\end{corollary}

It is clear that the functions%
\begin{equation}
\varphi(t)=\left\{
\begin{array}
[c]{c}%
0\text{ if }t\leq\tau,\\
u\left(  t\right)  \text{ if }t\geq\tau,
\end{array}
\right.  \label{Connection2}%
\end{equation}
where $u\left(  \text{\textperiodcentered}\right)  $ is the unique solution to
problem (\ref{LinearAut}) in $[\tau,+\infty)$ with $u_{\tau}=0$, are bounded
complete trajectories of $\mathcal{R}_{0}^{+}$ connecting $0$ and $v_{1}^{+}$,
as in this case the convergence in (\ref{ConvergBackwards}) holds true for
$t\rightarrow+\infty$. In fact, these are the unique possible connections.

\begin{theorem}
\label{StrAttr}Let (\ref{CondAttr}) hold. Any bounded complete trajectory of
$\mathcal{R}_{0}^{+}$ distinct from $0$ and $v_{1}^{+}$ is of the type
(\ref{Connection2}).
\end{theorem}

\begin{proof}
In view of Corollary \ref{Backward}, if $v\in\mathcal{A}^{+}$, then either
$v\equiv0$ or $v$ is positive. Let $\varphi\left(  \text{\textperiodcentered
}\right)  $ be a bounded complete trajectory of $\mathcal{R}_{0}^{+}$ distinct
from $0$ and $v_{1}^{+}$, so there is $t_{1}\in\mathbb{R}$ such that
$\varphi\left(  t_{1}\right)  $ is positive. By Lemmas
\ref{UniqueNondegenerate} and \ref{PositiveSolutions} there must be a first
time $\tau<t_{1}$ for which $\varphi\left(  \tau\right)  \equiv0$ and
$\gamma\left(  t\right)  >0$ for all $t>\tau$, and Corollary \ref{Backward2}
implies then that $\varphi\left(  t\right)  \equiv0$ for all $t\leq\tau$.
Hence, $\varphi$ is of the type (\ref{Connection2}).
\end{proof}

\bigskip

We obtain also that the trajectories inside the global attractor are regular.

\begin{proposition}
\label{RegularityAttr}Let (\ref{CondAttr}) hold. For any bounded complete
trajectory $\varphi$ of $\mathcal{R}_{0}^{+}$ which is not a fixed point, that
is, of the type (\ref{Connection2}), we have:

\begin{enumerate}
\item For any $t_{0}>\tau$ the function $u=\varphi\mid_{t\geq t_{0}}$ is the
unique non-negative solution to problem (\ref{Incl}) in $[t_{0},+\infty)$ with
$u_{t_{0}}=\varphi(t_{0});$

\item $\varphi\in C^{\infty}\left(  \left(  \tau,+\infty\right)  \times
\lbrack0,1]\right)  $ and $\varphi\in C((-\infty,+\infty),C^{1}([0,1]));$

\item There exists a time $t_{0}>\tau$ such that $u=\varphi\mid_{t\geq t_{0}}$
is the unique solution to problem (\ref{Incl}) in $[t_{0},+\infty)$ with
$u_{t_{0}}=\varphi(t_{0}).$
\end{enumerate}
\end{proposition}

\begin{proof}
The first two statements are a direct consequence of Lemmas
\ref{PositiveSolutions}, \ref{H3Regularity} and (\ref{ContDeriv}).

Taking into account that $V^{2r}\subset C^{1}([0,1])$ for $r>\frac{3}{4}$, we
obtain that
\[
\varphi(t)\rightarrow v_{1}^{+}\text{ in }C^{1}([0,1])\text{ as }%
t\rightarrow+\infty.
\]
Thus, there exist $t_{0}$ such that $\varphi_{x}(t,0)>0,\ \varphi
_{x}(t,1)<0,\ \varphi\left(  t,x\right)  >0$ if $x\in\left(  0,1\right)  $,
for all $t\geq t_{0}$. The result follows then from Lemma
\ref{PositiveSolutions3}.
\end{proof}

\begin{corollary}
\label{RegularityAttr2}Let (\ref{CondAttr}) hold. Then, $\mathcal{A}%
^{+}\subset C^{\infty}([0,1]).$
\end{corollary}

\subsection{The nonautonomous case\label{Nonautonomous}}

Throughout this section we assume that the functions $b\left(
\text{\textperiodcentered}\right)  ,\omega\left(  \text{\textperiodcentered
}\right)  $ satisfy (\ref{Cond1}) and, additionally, that
\begin{equation}
\omega_{1}<\pi^{2}. \label{CondAttr2}%
\end{equation}
Let $\mathcal{D}_{\tau}(u_{\tau})$ be the set of all solutions of problem
(\ref{Incl}) with initial condition $u_{\tau}$ at time $\tau$ and let
$\mathcal{R}_{\tau}=\cup_{u_{\tau}\in H}\mathcal{D}_{\tau}(u_{\tau})$,
$\mathcal{R=\cup}_{\tau\in\mathbb{R}}\mathcal{R}_{\tau}$. Denoting
$\mathbb{R}_{\geq}^{2}=\{\left(  t,s\right)  \in\mathbb{R}^{2}:t\geq s\}$ we
define the family of operators $U:\mathbb{R}_{\geq}^{2}\times H\rightarrow
P(H)$ given by%
\[
U(t,s,x)=\{u\left(  t\right)  :u\left(  \text{\textperiodcentered}\right)
\in\mathcal{D}_{s}(x)\}.
\]
$U$ is a strict multivalued process, that is, $U(t,t,$\textperiodcentered$)$
is the identity map and $U(t,s,U(s,\tau,x))=U(t,\tau,x)$ for all $t\geq
s\geq\tau$ and $x\in H$. Under assumptions (\ref{Cond1}), (\ref{CondAttr2})
this process has a compact strictly invariant pullback attractor
$\{\mathcal{A}(t)\}_{t\in\mathbb{R}}$ \cite[Theorem 5]{CVL19}, which means that:

\begin{itemize}
\item The sets $\mathcal{A}(t)$ are compact for any $t;$

\item $\mathcal{A}(t)=U(t,s,\mathcal{A}(s))$ for all $t\geq s$ (strict invariance);

\item $dist(U(t,s,B),\mathcal{A}(t))\rightarrow0$ as $s\rightarrow-\infty$ for
any bounded set $B\subset H$ (pullback attraction);

\item $\{\mathcal{A}(t)\}_{t\in\mathbb{R}}$ is the minimal pullback attracting
family, that is, if $\{\mathcal{K}(t)\}_{t\in\mathbb{R}}$ is a family of
closed pullback attracting sets, then $\mathcal{A}(t)\subset\mathcal{K}(t)$
for any $t.$
\end{itemize}

The function $\gamma:\mathbb{R}\rightarrow H$ is called a complete trajectory
of $\mathcal{R}$ if
\[
\varphi=\gamma\mid_{\lbrack\tau,+\infty)}\in\mathcal{R}_{\tau}\text{ for any
}\tau\in\mathbb{R}.
\]
A complete trajectory $\gamma$ is bounded if the set $\cup_{t\in\mathbb{R}%
}\gamma\left(  t\right)  $ is bounded in $H$. The pullback attractor
$\mathcal{A}(t)$ consists of the union of the elements of all the bounded
complete trajectories \cite[Lemma 6]{CVL19}, that is,%
\[
\mathcal{A}(t)=\{\gamma(t):\gamma\text{ is a bounded complete trajectory of
}\mathcal{R}\}.
\]
Moreover, it is known \cite[Theorem 5 and Corollary 6]{CVL19} that $\cup
_{t\in\mathbb{R}}\mathcal{A}(t)$ is bounded in $V$ (so $\overline{\cup
_{t\in\mathbb{R}}\mathcal{A}(t)}$ is compact in $H$) and that the sets
$\mathcal{A}(t)$ are compact in $V$. Let us prove that in fact the union
$\cup_{t\in\mathbb{R}}\mathcal{A}(t)$ is relatively compact in $V^{2r}$ for
$0\leq r<1.$

\begin{lemma}
\label{CompactV2rNonAut}Assume that (\ref{Cond1}), (\ref{CondAttr2}) hold.
Then the set $\cup_{t\in\mathbb{R}}\mathcal{A}(t)$ is relatively compact in
$V^{2r}$ for any $0\leq r<1$.
\end{lemma}

\begin{proof}
For any $z\in\mathcal{A}(t)$ the equality $\mathcal{A}(t)=U(t,t-1,\mathcal{A}%
(t-1))$ implies the existence of $u\left(  \text{\textperiodcentered}\right)
\in\mathcal{D}_{t-1}(u_{t-1})$, $u_{t-1}\in\mathcal{A}(t-1)$, such that
$z=u(t)$ and that $u\left(  s\right)  \in\mathcal{A}(s)$ for any $s\geq t-1.$
As $u$ is a mild solution, by the variation of constants formula we get%
\[
z=e^{-A}u_{t-1}+\int_{t-1}^{t}e^{-A(t-s)}f(s)ds,
\]
where $f(s,x)\in b(t)H_{0}(u(t,x))+\omega(t)u(t,x)$ for a.a. $(t,x)$. By the
boundedness of $\cup_{s\in\mathbb{R}}\mathcal{A}(s)$ in $H,$ there exists an
universal constant $C>0$ such that
\[
\left\Vert u_{t-1}\right\Vert \leq C,\ \left\Vert f\right\Vert _{L^{\infty
}(t-1,t;H)}\leq C,
\]
where $C$ does not depend either on the chosen $z$ or $t$. Then using the
estimates of the norm of $e^{-At}$ in $V^{2r}$ \cite[Theorem 37.5]{SellBook}
there exist $M_{r}>0$, $a\in\mathbb{R}$ such that%
\begin{align*}
\left\Vert A^{r}z\right\Vert  &  \leq\left\Vert A^{r}e^{-A}u_{t-1}\right\Vert
+\int_{t-1}^{t}\left\Vert A^{r}e^{-A(t-s)}f(s)\right\Vert ds\\
&  \leq M_{r}e^{-a}C+M_{r}C\int_{t-1}^{t}\left(  t-s\right)  ^{-r}ds\leq
C_{r}.
\end{align*}
Therefore, $\mathcal{A}(t)$ is bounded in $V^{2r}$ if $r<1$ uniformly in
$t\in\mathbb{R}$. The compact embedding $V^{\alpha}\subset V^{\beta}$ for
$\alpha>\beta$ implies that $\cup_{t\in\mathbb{R}}\mathcal{A}(t)$ is
relatively compact in $V^{2r}$ for any $0\leq r<1$.
\end{proof}

\bigskip

By \cite[Theorem 6]{CVL19} we know that under assumptions (\ref{Cond1}),
(\ref{CondAttr2}) there exists a bounded complete trajectory $\xi_{M}$ such that:

\begin{enumerate}
\item $\xi_{M}(t)$ is positive for any $t\in\mathbb{R}$.

\item $-\xi_{M}(t)\leq\gamma\left(  t\right)  \leq\xi_{M}(t)$ for all
$t\in\mathbb{R}$ and any bounded complete trajectory $\gamma.$

\item $w_{b_{0},\omega_{0}}^{+}\leq\xi_{M}(t)\leq w_{b_{1},\omega_{1}}^{+}$,
where $w_{b_{i},\omega_{i}}^{+}$ denote the positive fixed point $v_{1}^{+}$
of the autonomous problem (\ref{Incl}) with $b=b_{i},\ \omega=\omega_{i}$.

\item $\xi_{M}$ is the unique bounded complete trajectory such that $\xi
_{M}(t)\in D$ for all $t\in\mathbb{R}$, where we recall that $D=\{v\in
V:v(x)>0$ for $x\in\left(  0,1\right)  \}$.
\end{enumerate}

From these properties the following lemma follows immediately.

\begin{lemma}
\label{Backward4}Assume that (\ref{Cond1}), (\ref{CondAttr2}) hold. If
$\gamma$ is a bounded complete trajectory such that $\gamma\left(  t\right)
\rightarrow0$ as $t\rightarrow-\infty$, then there exists a time $t_{0}$ such
that $\gamma\left(  t_{0}\right)  \not \in D.$
\end{lemma}

\begin{proof}
By contradiction let $\gamma$ be a bounded complete trajectory such that
$\gamma\left(  t\right)  \rightarrow0$ as $t\rightarrow-\infty$ and
$\gamma\left(  t\right)  \in D$ for any $t\in\mathbb{R}.$ Hence, $\gamma
\equiv\xi_{M}$. However, $\xi_{M}(t)\geq w_{b_{0},\omega_{0}}^{+}$, so
$\gamma(t)$ cannot converge to $0$.
\end{proof}

\bigskip

Now we are in position to improve the fourth point above and obtain that in
fact $\xi_{M}$ is the unique non-degenerate bounded complete solution at
$-\infty$, which means that there exists a time $t_{0}$ such that $\xi
_{M}(t)\in D$ for any $t\leq t_{0}$. This question was left as an open problem
in \cite[Remark 2]{CVL19}.

\begin{lemma}
\label{NonautEquilibrium}Assume that (\ref{Cond1}), (\ref{CondAttr2}) hold. Then:

\begin{itemize}
\item If $\gamma\left(  \text{\textperiodcentered}\right)  $ is a non-negative
bounded complete trajectory such that for some $t_{0}$ we have that
$\gamma\left(  t\right)  \in D$ for all $t\leq t_{0}$, then $\gamma\left(
t\right)  =\xi_{M}(t)$ for all $t\in\mathbb{R}$, that is, $\xi_{M}$ is the
unique non-negative bounded complete trajectory which is non-degenerate at
$-\infty.$ Also, for any $\tau\in\mathbb{R}$ the function $u=\xi_{M}%
\mid_{t\geq\tau}$ is the unique non-negative solution to problem (\ref{Incl})
with $u_{\tau}=\xi_{M}(\tau)$.

\item Assume, additionally, that $b,\omega\in W_{loc}^{1,2}(\mathbb{R})$. Then
$\xi_{M}$ is the unique bounded complete trajectory which is non-degenerate at
$-\infty.$ Also, for any $\tau\in\mathbb{R}$ the function $u=\xi_{M}%
\mid_{t\geq\tau}$ is the unique solution to problem (\ref{Incl}) with
$u_{\tau}=\xi_{M}(\tau)$.
\end{itemize}
\end{lemma}

\begin{proof}
Let us consider the first statement. Since $\gamma$ is non-negative, Lemma
\ref{PositiveSolutions} implies that for all $\tau\in\mathbb{R}$ the function
$u=\gamma\mid_{t\geq\tau}$ is the unique non-negative solution to problem
(\ref{Incl}) and $\gamma\left(  \tau\right)  \in D$. Thus, by the previous
results it follows that $\gamma\left(  t\right)  =\xi_{M}(t)$ for all
$t\in\mathbb{R}$.

We prove further the second statement. As $\mathcal{A}(t)\subset C^{1}([0,1])$
by Lemma \ref{CompactV2rNonAut} and the embedding $V^{2r}\subset C^{1}([0,1])$
for $r>\frac{3}{4}$, it follows from $w_{b_{0},\omega_{0}}^{+}\leq\xi_{M}(t),$
$\frac{d}{dx}w_{b_{0},\omega_{0}}^{+}(0)>0$, $\frac{d}{dx}w_{b_{0},\omega_{0}%
}^{+}(1)<0$ that $\frac{\partial}{\partial x}\xi_{M}(t,0)>0,\ \frac{\partial
}{\partial x}\xi_{M}(t,1)<0$ for any $t\in\mathbb{R}$. Hence, from Lemma
\ref{PositiveSolutions3} and $w_{b_{0},\omega_{0}}^{+}\in D$ we have that for
any $\tau\in\mathbb{R}$ the function $u=\xi_{M}\mid_{t\geq\tau}$ is the unique
solution to problem (\ref{Incl}) with $u_{\tau}=\xi_{M}(\tau)$. Let now
$\gamma$ be a bounded complete trajectory which is non-degenerate at
$-\infty.$ Thus, by Corollary 8 in \cite{CVL19} there exists a time $t_{0}$
such that $\gamma\left(  t\right)  =\xi_{M}(t)$ for all $t\leq t_{0}$. By the
uniqueness of the solution $u=\xi_{M}\mid_{t\geq t_{0}}$ with $u_{t_{0}}%
=\xi_{M}(t_{0})$ we get that $\gamma\left(  t\right)  =\xi_{M}(t)$ for any
$t\in\mathbb{R}.$
\end{proof}

\bigskip

\begin{corollary}
Assume that (\ref{Cond1}), (\ref{CondAttr2}) hold. Then if $\gamma$ is a
non-negative bounded complete trajectory such that $\gamma\left(  t\right)
\rightarrow0$ as $t\rightarrow-\infty$, for any $t_{1}$ there exists
$t_{0}\leq t_{1}$ such that $\gamma\left(  t_{0}\right)  \not \in D.$

If, additionally, $b,\omega\in W_{loc}^{1,2}(\mathbb{R})$, then the result is
valid for any bounded complete trajectory $\gamma$ such that $\gamma\left(
t\right)  \rightarrow0$ as $t\rightarrow-\infty$.
\end{corollary}

\begin{proof}
If there is $t_{1}$ such that $\gamma\left(  t\right)  \in D$ for all $t\leq
t_{1}$, then Lemma \ref{NonautEquilibrium} implies that $\gamma\equiv\xi_{M},$
so that the convergence $\gamma\left(  t\right)  \rightarrow0$, as
$t\rightarrow-\infty$, is not possible. The same proof is valid for the second statement.
\end{proof}

\bigskip

We denote by $\mathcal{D}_{\tau}^{+}(u_{\tau})$ the set of all non-negative
solutions of problem (\ref{Incl}) with initial condition $u_{\tau}\in H^{+}$
at time $\tau$ and let $\mathcal{R}_{\tau}^{+}=\cup_{u_{\tau}\in H}%
\mathcal{D}_{\tau}^{+}(u_{\tau})$, $\mathcal{R}^{+}\mathcal{=\cup}_{\tau
\in\mathbb{R}}\mathcal{R}_{\tau}^{+}$. We define the map $U^{+}:\mathbb{R}%
_{\geq}\times H^{+}\rightarrow P(H^{+})$ given by%
\[
U^{+}(t,\tau,u_{\tau})=\{u(t):u\in\mathcal{D}_{\tau}^{+}(u_{\tau})\}.
\]
$U^{+}$ is a strict multivalued process, which follows easily by the
translation and concatenation properties of solutions (see \cite[Section
4]{CVL19}). The bounded complete trajectories of $\mathcal{R}^{+}$ are all the
bounded complete trajectories $\gamma\left(  \text{\textperiodcentered
}\right)  $ of $\mathcal{R}$ such that $\gamma(t)\geq0$ for any $t\in
\mathbb{R}$. The process $U^{+}$ possesses a compact strictly invariant
pullback attractor, denoted by $\{\mathcal{A}^{+}(t)\}_{t\in\mathbb{R}}$, and%
\[
\mathcal{A}^{+}(t)=\{\gamma(t):\gamma\text{ is a bounded complete trajectory
of }\mathcal{R}^{+}\}.
\]
As $\mathcal{A}^{+}(t)\subset\mathcal{A}(t)$ for all $t$, $\{\mathcal{A}%
^{+}(t)\}_{t\in\mathbb{R}}$ inherits the regularity properties of
$\{\mathcal{A}(t)\}_{t\in\mathbb{R}}$ described above.

We have seen in Lemma \ref{NonautEquilibrium} that the function $\xi_{M}$ is
the unique non-negative bounded complete trajectory which is non-degenerate at
$-\infty$. This function plays the role of a positive nonautonomous
equilibrium. As in the autonomous case, we will obtain a precise
characterization of the pullback attractor in the positive cone in terms of
the bounded complete trajectories connecting $0$ and $\xi_{M}$.

As in (\ref{Connection2}), we define the map
\begin{equation}
\varphi(t)=\left\{
\begin{array}
[c]{c}%
0\text{ if }t\leq\tau,\\
u\left(  t\right)  \text{ if }t\geq\tau,
\end{array}
\right.  \label{Connection3}%
\end{equation}
where $u\left(  \text{\textperiodcentered}\right)  $ is the unique solution to
problem
\begin{equation}
\left\{
\begin{array}
[c]{l}%
\dfrac{\partial u}{\partial t}-\dfrac{\partial^{2}u}{\partial x^{2}}\in
b(t)+\omega(t)u,\text{ on }(\tau,+\infty)\times(0,1),\\
u|_{\partial\Omega}=0,\\
u(\tau,x)=0,
\end{array}
\right.  \label{LinearNonAut}%
\end{equation}
It is clear that a function of the type (\ref{Connection3}) is a bounded
complete trajectory of $\mathcal{R}^{+}$. Alos, it is easy to see that
$\left\Vert u\left(  t\right)  -\xi_{M}(t)\right\Vert \rightarrow0$ as
$t\rightarrow+\infty$. Indeed, for $v(t)=u(t)-\xi_{M}(t)$ we have that%
\[
\frac{d}{dt}\left\Vert v\right\Vert ^{2}+2\pi^{2}\left\Vert v\right\Vert
^{2}\leq2\omega\left(  t\right)  \left\Vert v\right\Vert ^{2}\leq2\omega
_{1}\left\Vert v\right\Vert ^{2}.
\]
We denote $\delta=\pi^{2}-\omega_{1}>0$. Then%
\[
\left\Vert v(t)\right\Vert ^{2}\leq e^{-2\delta(t-\tau)}\left\Vert v\left(
\tau\right)  \right\Vert ^{2}\rightarrow0\text{ as }t\rightarrow+\infty.
\]
Obviously, $u\left(  t\right)  \rightarrow0$ as $t\rightarrow-\infty$. Hence,
$\varphi\left(  \text{\textperiodcentered}\right)  $ is a connection from $0$
to $\xi_{M}$, that is,%
\begin{align*}
u\left(  t\right)   &  \rightarrow0\text{ as }t\rightarrow-\infty,\\
u\left(  t\right)  -\xi_{M}(t)  &  \rightarrow0\text{ as }t\rightarrow+\infty.
\end{align*}
We will establish that these are the only possible bounded complete
trajectories of $\mathcal{R}^{+}$.

\begin{theorem}
\label{StrPullbackAttr}Assume that (\ref{Cond1}), (\ref{CondAttr2}) hold. Then
any bounded complete trajectory of $\mathcal{R}^{+}$ distinct from $0$ and
$\xi_{M}$ is of the type (\ref{Connection3}).
\end{theorem}

\begin{proof}
By Corollary \ref{Backward}, if $v\in\mathcal{A}^{+}(t)$, then either
$v\equiv0$ or $v$ is positive. We take a bounded complete trajectory
$\gamma\left(  \text{\textperiodcentered}\right)  $ of $\mathcal{R}^{+}$
different from $0$ and $\xi_{M}$, so there exists $t_{1}\in\mathbb{R}$ such
that $\gamma\left(  t_{1}\right)  >0$. By Lemmas \ref{PositiveSolutions} and
\ref{NonautEquilibrium} there exists a first time $\tau<t_{1}$ for which
$\gamma\left(  \tau\right)  \equiv0$ and $\gamma\left(  t\right)  >0$ for all
$t>\tau$. Thus, Corollary \ref{Backward2} implies that $\gamma\left(
t\right)  \equiv0$ for all $t\leq\tau$, and then $\gamma$ is of the type
(\ref{Connection3}).
\end{proof}

\begin{remark}
The structure of the pullback attractor in the positive cone is the same as
for the ordinary differential inclusion
\[
\left\{
\begin{array}
[c]{l}%
\dfrac{du}{dt}+\lambda u\in b(t)H_{0}(u),\text{ on }(\tau,+\infty),\\
u(\tau)=u_{\tau},
\end{array}
\right.
\]
where $\lambda>0,$ which was studied in \cite{CVL16}.
\end{remark}

\bigskip

We conclude this section by obtaining some regularity results for the
solutions inside the pullback attractor. The following proposition is proved
in the same way as Proposition \ref{RegularityAttr} but using that any bounded
complete trajectory $\gamma$ of the type (\ref{Connection3}) satisfies that%
\[
\left\Vert \gamma\left(  t\right)  -\xi_{M}(t)\right\Vert _{C^{1}%
([0,1])}\rightarrow0\text{ as }t\rightarrow+\infty.
\]

\begin{proposition}
\label{RegularityAttr3}Assume that (\ref{Cond1}), (\ref{CondAttr2}) hold. Then
any bounded complete trajectory $\gamma$ of $\mathcal{R}^{+}$ different from
$0$ and $\xi_{M}$, that is, of the type (\ref{Connection3}), satisfies that
for any $t_{0}>\tau$ the function $u=\gamma\mid_{t\geq t_{0}}$ is the unique
non-negative solution to problem (\ref{Incl}) in $[t_{0},+\infty)$ with
$u_{t_{0}}=\gamma(t_{0}).$ If, additionally, $b,\omega\in W_{loc}%
^{k,2}(\mathbb{R})$, $k\in\mathbb{N}$, then:

\begin{enumerate}
\item $\gamma\in C^{k}((\tau,+\infty),V)$, $\gamma\in\cap_{j=0}^{k-1}%
C^{j}\left(  (\tau,+\infty),H^{2(k-j)+1}(0,1)\right)  $ and $\gamma\in
C^{k}(\left(  \tau,+\infty\right)  \times\lbrack0,1]).$

\item There exists a time $t_{0}>\tau$ such that $u=\gamma\mid_{t\geq t_{0}}$
is the unique solution to problem (\ref{Incl}) in $[t_{0},+\infty)$ with
$u_{t_{0}}=\gamma(t_{0}).$
\end{enumerate}

Moreover, if $b,\omega\in W_{loc}^{k,2}(\mathbb{R})$ for any $k\in\mathbb{N}$,
then $\gamma\in C^{\infty}(\left(  \tau,+\infty\right)  \times\lbrack0,1]).$
\end{proposition}

In the case of the nonautonomous equilibrium $\xi_{M},$ the regularity is
extended to the whole line $\left(  -\infty,+\infty\right)  .$

\begin{proposition}
Assume that (\ref{Cond1}), (\ref{CondAttr2}) hold. If $b,\omega\in
W_{loc}^{k,2}(\mathbb{R})$, $k\in\mathbb{N}$, then%
\[
\xi_{M}\in C^{k}((-\infty,+\infty),V),\ \xi_{M}\in\cap_{j=0}^{k-1}C^{j}\left(
(-\infty,+\infty),H^{2(k-j)+1}(0,1)\right)  ,\ \xi_{M}\in C^{k}(\left(
-\infty,+\infty\right)  \times\lbrack0,1]).
\]

Moreover, if $b,\omega\in W_{loc}^{k,2}(\mathbb{R})$ for any $k\in\mathbb{N}$,
then $\xi_{M}\in C^{\infty}(\left(  -\infty,+\infty\right)  \times
\lbrack0,1]).$
\end{proposition}

\begin{corollary}
\label{RegularityAttr4}Assume that (\ref{Cond1}), (\ref{CondAttr2}) hold. We have:

\begin{enumerate}
\item If $b,\omega\in W_{loc}^{k,2}(\mathbb{R}),$ where $k\in\mathbb{N}$, then
$\cup_{t\in\mathbb{R}}\mathcal{A}^{+}(t)\subset H^{2k+1}(0,1)\subset
C^{2k}([0,1]).$

\item If $b,\omega\in W_{loc}^{k,2}(\mathbb{R})$ for any $k\in\mathbb{N}$,
then $\cup_{t\in\mathbb{R}}\mathcal{A}^{+}(t)\subset C^{\infty}([0,1]).$
\end{enumerate}
\end{corollary}

\section{Comparison of different non-autonomous attractors}

In this section we are going to consider different type of attractors like the
uniform attractor, the cocycle attractor, the skew-product semiflow attractor
and, of course, the pullback attractor in order to establish the relationship
between them. Moreover, we will prove that the attractor of the skew-product
semiflow in the positive cone has a gradient structure.

For this aim we need to consider not only problem (\ref{Incl}) itself but all
the problems generated by the translations of the functions $b\left(
\text{\textperiodcentered}\right)  ,\ \omega\left(  \text{\textperiodcentered
}\right)  $, that is, by the hull of these functions. In order to define this
hull properly, along with conditions (\ref{Cond1}), (\ref{CondAttr2}), we need
the following extra assumption:%
\begin{equation}
\text{the functions }b\text{ and }\omega\text{ are uniformly continuous in
}\left(  -\infty,+\infty\right)  . \label{CondAttr3}%
\end{equation}

For a sequence $t_{n}\in\mathbb{R}$ we define the functions $\sigma_{n}\left(
\text{\textperiodcentered}\right)  =\left(  b\left(  \text{\textperiodcentered
}+t_{n}\right)  ,\omega\left(  \text{\textperiodcentered}+t_{n}\right)
\right)  \in C(\mathbb{R},\mathbb{R}^{2})$. As usual, the space $C(\mathbb{R}%
,\mathbb{R}^{2})$ is equipped with the topology of uniform convergence on
compact sets of $\mathbb{R}$. The sequence $\sigma_{n}$ is uniformly bounded
and uniformly (in $\left(  -\infty,+\infty\right)  $) equicontinuous.
Therefore, by applying the Ascoli-Arzel\`{a} theorem and a diagonal argument
there exists a subsequence $\sigma_{n^{\prime}}$ and a uniformly continuous
function $\sigma\left(  \text{\textperiodcentered}\right)  =\left(
\overline{b}\left(  \text{\textperiodcentered}\right)  ,\overline{\omega
}\left(  \text{\textperiodcentered}\right)  \right)  $ satisfying
(\ref{Cond1}), (\ref{CondAttr2}) as well (with the same constants
$b_{0},\ b_{1},\ \omega_{1}$) such that
\[
\sigma_{n}\rightarrow\sigma\text{ in }C(\mathbb{R},\mathbb{R}^{2})\text{.}%
\]
Thus, we define the hull of $\left(  b\left(  \text{\textperiodcentered
}\right)  ,\omega\left(  \text{\textperiodcentered}\right)  \right)  $ by
\[
\Sigma=cl_{C(\mathbb{R},\mathbb{R}^{2})}\{\left(  b\left(
\text{\textperiodcentered}+t\right)  ,\omega\left(  \text{\textperiodcentered
}+t\right)  \right)  :t\in\mathbb{R}\},
\]
which is a compact set of the metrizable space $C(\mathbb{R},\mathbb{R}^{2}),$
and consider the family of problems (\ref{Incl}) given by each element
$\sigma\left(  \text{\textperiodcentered}\right)  =\left(  \overline{b}\left(
\text{\textperiodcentered}\right)  ,\overline{\omega}\left(
\text{\textperiodcentered}\right)  \right)  \in\Sigma$. We define the
translation operator $\theta_{s}:\Sigma\rightarrow\Sigma,\ s\in\mathbb{R},$
given by $\theta_{s}\sigma\left(  \text{\textperiodcentered}\right)
=\sigma\left(  \text{\textperiodcentered}+s\right)  $. Clearly, $\theta_{0}$
is the identity operator and $\theta_{s+r}=\theta_{s}\circ\theta_{r}$, so
$\theta$ is a group, which is called the driving group in $\Sigma$. Also, it
is not difficult to see that the map $\left(  t,\sigma\right)  \mapsto
\theta_{t}\sigma$ is continuous and that $\theta_{t}\Sigma=\Sigma$ for any
$t\in\mathbb{R}$.

Let now $\mathcal{D}_{\tau,\sigma}(u_{\tau})$ be the set of all solutions of
problem (\ref{Incl}) with initial condition $u_{\tau}$ at time $\tau$ with
symbol $\sigma$ and let $\mathcal{R}_{\tau,\sigma}=\cup_{u_{\tau}\in
H}\mathcal{D}_{\tau,\sigma}(u_{\tau})$, $\mathcal{R}_{\sigma}=\mathcal{\cup
}_{\tau\in\mathbb{R}}\mathcal{R}_{\tau,\sigma}$. For each $\sigma\in\Sigma$
the operator $U_{\sigma}:\mathbb{R}_{\geq}^{2}\times H\rightarrow P(H)$ given
by%
\[
U_{\sigma}(t,\tau,x)=\{u\left(  t\right)  :u\left(  \text{\textperiodcentered
}\right)  \in\mathcal{D}_{\tau,\sigma}(x)\}
\]
is a strict multivalued process. We obtain then a family of multivalued
processes. Additionally, we will need the following translation property:%
\begin{equation}
U_{\sigma}(t+h,\tau+h,x)=U_{\theta_{h}\sigma}(t,\tau,x)\text{ for all }x\in
H\text{, }\sigma\in\Sigma\text{, }h\in\mathbb{R}\text{, }\left(
t,\tau\right)  \in\mathbb{R}_{\geq}^{2}. \label{Translation}%
\end{equation}

\begin{lemma}
\label{TranslationLemma}Let (\ref{Cond1}), (\ref{CondAttr2}), (\ref{CondAttr3}%
) hold. Then the family $\{U_{\sigma}\}$ satisfies (\ref{Translation}).
\end{lemma}

\begin{proof}
If $y\in U_{\sigma}(t+h,\tau+h,x)$, there is $u\in\mathcal{D}_{\tau+h,\sigma
}(x)$ such that $y=u(t+h)$, where $\sigma=\left(  \overline{b},\overline
{\omega}\right)  $. Using the concept of mild solution we have a selection
$r\in L_{loc}^{2}(\tau+h,+\infty;H)$ of the map $R\left(  t,u\left(  t\right)
\right)  $ such that%
\[
u(t+h)=e^{-A(t-\tau)}x+\int_{\tau+h}^{t+h}e^{-A(t+h-s)}r(s)ds\text{ for }%
\tau\leq t<+\infty.
\]
Let $v\left(  \text{\textperiodcentered}\right)  =u\left(
\text{\textperiodcentered}+h\right)  ,\ r_{h}\left(  \text{\textperiodcentered
}\right)  =r\left(  \text{\textperiodcentered}+h\right)  =\theta_{h}r\in
L_{loc}^{2}(\tau,+\infty;H)$. Then%
\begin{align*}
v\left(  t\right)   &  =e^{-A(t-\tau)}x+\int_{\tau}^{t}e^{-A(t-s)}r(s+h)ds\\
&  =e^{-A(t-\tau)}x+\int_{\tau}^{t}e^{-A(t-s)}r_{h}(s)ds\text{ for all }%
\tau\leq t<+\infty.
\end{align*}
As $r_{h}(t,x)\in\overline{b}(t+h)H_{0}\left(  v\left(  t,x\right)  \right)
+\overline{\omega}(t+h)v\left(  t,x\right)  $ for a.a. $\left(  t,x\right)  $,
we infer that $v\left(  \text{\textperiodcentered}\right)  $ is a mild
solution on $\left(  \tau,+\infty\right)  $ for the symbol $\theta_{h}\sigma$,
so that $v\in\mathcal{D}_{\tau,\theta_{h}\sigma}(x)$. Hence, $y\in
U_{\theta_{h}\sigma}(t,\tau,x)$.

The converse inclusion is proved in a similar way.
\end{proof}

\bigskip

By the results in the previous section we known that each $U_{\sigma}$
possesses a compact strictly invariant pullback attractor $\{\mathcal{A}%
_{\sigma}(t)\}_{t\in\mathbb{R}}$, which can be characterized in terms of the
bounded complete trajectories of $\mathcal{R}_{\sigma}:$%
\[
\mathcal{A}_{\sigma}(t)=\{\gamma(t):\gamma\text{ is a bounded complete
trajectory of }\mathcal{R}_{\sigma}\}.
\]

We can extend Lemma \ref{CompactV2rNonAut} by taking the union of all the
attractors over $\Sigma$.

\begin{lemma}
\label{CompactV2rNonAut2}Assume that (\ref{Cond1}), (\ref{CondAttr2}),
(\ref{CondAttr3}) hold. Then the set $\cup_{\sigma\in\Sigma,t\in\mathbb{R}%
}\mathcal{A}_{\sigma}(t)$ is relatively compact in $V^{2r}$ for any $0\leq
r<1$.
\end{lemma}

\begin{proof}
The proof is the same as in Lemma \ref{CompactV2rNonAut} by taking into
account that the\ constant $C$ in the estimate $\left\Vert f\right\Vert
_{L^{\infty}(t-1,t;H)}\leq C$ is independent of $\sigma\in\Sigma.$
\end{proof}

\bigskip

As before, we extend the above results when we restrict ourselves to
non-negative solutions. Let $\mathcal{D}_{\tau,\sigma}^{+}(u_{\tau})$ be the
set of all non-negative solutions of problem (\ref{Incl}) with initial
condition $u_{\tau}$ at time $\tau$ with symbol $\sigma$ and let
$\mathcal{R}_{\tau,\sigma}^{+}=\cup_{u_{\tau}\in H}\mathcal{D}_{\tau,\sigma
}^{+}(u_{\tau})$, $\mathcal{R}_{\sigma}^{+}=\mathcal{\cup}_{\tau\in\mathbb{R}%
}\mathcal{R}_{\tau,\sigma}$. For each $\sigma\in\Sigma$ the operator
$U_{\sigma}^{+}:\mathbb{R}_{\geq}^{2}\times H^{+}\rightarrow P(H^{+})$ is
given by%
\[
U_{\sigma}^{+}(t,\tau,x)=\{u\left(  t\right)  :u\left(
\text{\textperiodcentered}\right)  \in\mathcal{D}_{\tau,\sigma}^{+}(x)\}.
\]
$U_{\sigma}^{+}$ is a strict multivalued process and satisfies the translation
property (\ref{Translation}). The process $U_{\sigma}^{+}$ possesses a compact
strictly invariant pullback attractor $\{\mathcal{A}_{\sigma}^{+}%
(t)\}_{t\in\mathbb{R}}$ and%
\[
\mathcal{A}_{\sigma}^{+}(t)=\{\gamma(t):\gamma\text{ is a bounded complete
trajectory of }\mathcal{R}_{\sigma}^{+}\}.
\]
As proved in the previous section, for each $\sigma$ the attractor
$\{\mathcal{A}_{\sigma}^{+}(t)\}$ consists of the equilibria $0$ and
$\xi_{M,\sigma}\left(  t\right)  $ and the bounded complete trajectories,
given by (\ref{Connection3}), which connect them.

\subsection{The cocycle attractor}

We define for the family $\{U_{\sigma}\}$ the map $\phi:\mathbb{R}^{+}%
\times\Sigma\times H\rightarrow P\left(  H\right)  $ given by%
\[
\phi(t,\sigma,x)=U_{\sigma}(t,0,x)\text{ for all }t\geq0\text{, }\sigma
\in\Sigma,\ x\in H,
\]
which is a strict multivalued cocycle, that is, $\phi\left(  0,\sigma
,\text{\textperiodcentered}\right)  $ is the identity map and we have that
$\phi\left(  t+s,\sigma,x\right)  =\phi\left(  t,\theta_{s}\sigma
,\phi(s,\sigma,x\right)  )$ for all $t\geq s\geq0$, $\sigma\in\Sigma$, $x\in
H$ \cite[Proposition 1]{Cui}.

By \cite[Corollary 4]{Cui} the family $\{\mathcal{A}(\sigma)\}_{\sigma
\in\Sigma}$ defined by
\[
\mathcal{A}(\sigma)=\mathcal{A}_{\sigma}(0),
\]
where $\{\mathcal{A}_{\sigma}(t)\}_{t\in\mathbb{R}}$ is the pullback attractor
of $U_{\sigma}$, is a compact strictly invariant cocycle attractor for $\phi$,
which means that:

\begin{itemize}
\item The sets $\mathcal{A}(\sigma)$ are compact for any $\sigma\in\Sigma.$

\item $\mathcal{A}\left(  \theta_{t}\sigma\right)  =\phi\left(  t,\sigma
,\mathcal{A}\left(  \sigma\right)  \right)  $ for all $\sigma\in\Sigma$,
$t\geq0$ (strict invariance).

\item $\lim_{t\rightarrow+\infty}dist(\phi\left(  t,\theta_{-t}\sigma
,B\right)  ,\mathcal{A}(\sigma))=0$ for any bounded set $B\subset H$ and all
$\sigma\in\Sigma$ (pullback attraction).

\item $\{\mathcal{A}(\sigma)\}_{\sigma\in\Sigma}$ is minimal, that is, if
$\{\mathcal{A}^{\prime}(\sigma)\}_{\sigma\in\Sigma}$ is a family of closed
sets satisfying the pullback attraction property, then $\mathcal{A}%
(\sigma)\subset\mathcal{A}^{\prime}(\sigma)$ for any $\sigma\in\Sigma.$
\end{itemize}

Moreover,
\[
\mathcal{A}(\theta_{t}\sigma)=\mathcal{A}_{\sigma}(t)\text{ for all }\sigma
\in\Sigma\text{, }t\in\mathbb{R},
\]
so%
\[
\mathcal{A}(\theta_{t}\sigma)=\{\gamma\left(  t\right)  :\gamma\text{ is
bounded complete trajectory of }\mathcal{R}_{\sigma}\}.
\]
Also, Lemma \ref{CompactV2rNonAut2} implies that $\cup_{\sigma\in\Sigma
}\mathcal{A}(\sigma)$ is relatively compact in $V^{2r}$ for any $0\leq r<1.$

In the same way, we define the cocycle for non-negative solutions $\phi
^{+}:\mathbb{R}^{+}\times\Sigma\times H^{+}\rightarrow P\left(  H^{+}\right)
$ by
\[
\phi^{+}(t,\sigma,x)=U_{\sigma}^{+}(t,0,x)\text{ for all }t\geq0\text{,
}\sigma\in\Sigma,\ x\in H,
\]
which possesses a compact strictly invariant cocycle attractor $\{\mathcal{A}%
_{\sigma}^{+}(t)\}_{t\in\mathbb{R}}$ which satisfies%
\[
\mathcal{A}^{+}(\theta_{t}\sigma)=\mathcal{A}_{\sigma}^{+}(t)\text{ for all
}\sigma\in\Sigma\text{, }t\in\mathbb{R},
\]
and%
\[
\mathcal{A}^{+}(\theta_{t}\sigma)=\{\gamma\left(  t\right)  :\gamma\text{ is a
bounded non-negative complete trajectory of }\mathcal{R}_{\sigma}\}.
\]

We finish this subsection by noticing (see the proof of Lemma
\ref{TranslationLemma}) that $\gamma$ is a complete trajectory of
$\mathcal{R}_{\sigma}$ if and only if
\[
\gamma\left(  \text{\textperiodcentered}+s\right)  \mid_{\lbrack0,+\infty)}%
\in\mathcal{R}_{0,\theta_{s}\sigma}\text{ for any }s\in\mathbb{R}.
\]

\subsection{The skew-product semiflow attractor}

We will denote by $\mathcal{X}$ the product space $H\times\Sigma$ with the
metric $\rho_{\mathcal{X}}$ given by
\[
\rho_{\mathcal{X}}\left(  \left(  x_{1},\sigma_{1}\right)  ,\left(
x_{2},\sigma_{2}\right)  \right)  =\left\Vert x_{1}-x_{2}\right\Vert
+\rho\left(  \sigma_{1},\sigma_{2}\right)  ,
\]
where $\rho$ is a metric in the space $C\left(  \mathbb{R},\mathbb{R}%
^{2}\right)  $. Also, let $\mathcal{P}_{H}:\mathcal{X}\rightarrow H$ be the
projector onto $H$, that is, for a subset $C\subset\mathcal{X}$ we put
\[
\mathcal{P}_{H}(C)=\{u\in H:\left(  u,\sigma\right)  \in C\text{ for some
}\sigma\in\Sigma\}.
\]

From the cocycle $\phi$ we define the skew product semiflow $\Pi
:\mathbb{R}^{+}\times\mathcal{X}\rightarrow P(\mathcal{X})$ given by
\[
\Pi(t,\left(  x,\sigma\right)  )=\left(  \phi\left(  t,\sigma,x\right)
,\theta_{t}\sigma\right)  .
\]
Since $\phi$ is a strict cocycle it is easy to check that $\Pi$ is a strict
multivalued semiflow, which means that $\Pi\left(  0,\text{\textperiodcentered
}\right)  $ is the identity map and $\Pi(t,y)=\Pi(t,\Pi(s,y))$ for any
$y\in\mathcal{X}$, $t,s\geq0$.

We need to prove some properties of $\Pi$ leading to the existence of a global
attractor. After that we will establish the relationship with the cocycle
attractor and study its structure.

\begin{lemma}
\label{CompactAbsorbing}Let (\ref{Cond1}), (\ref{CondAttr2}), (\ref{CondAttr3}%
) hold. Then $\Pi$ possesses a compact absorbing set $\mathbb{K}$, which means
that for any bounded set $\mathbb{B}$ there exists $T(\mathbb{B})$ such that
$\Pi(t,\mathbb{B})\subset\mathbb{K}$ if $t\geq T.$
\end{lemma}

\begin{proof}
In a standard way (see \cite[Lemma 5]{CVL19}) for any solution to problem
(\ref{Incl}) in $[0,+\infty)$ we obtain the estimates%
\begin{equation}
\left\Vert u\left(  t\right)  \right\Vert ^{2}\leq e^{-\delta t}\left\Vert
u\left(  0\right)  \right\Vert ^{2}+\frac{C_{1}}{\delta}\text{ for all }%
t\geq0, \label{IneqL2}%
\end{equation}%
\begin{equation}
\int_{t-\alpha}^{t}\left\Vert \dfrac{\partial u}{\partial x}\right\Vert
^{2}dr\leq\frac{\pi^{2}C_{1}}{\delta}+\frac{\pi^{2}}{\delta}\left\Vert
u(t-\alpha)\right\Vert ^{2}, \label{IneqIntH1}%
\end{equation}
where $\alpha\in(0,1]$ is arbitrary and $C,\delta>0$ are universal constants
which are independent of $\sigma\in\Sigma$ (they depend only on the constants
$b_{1},\omega_{1}$ from (\ref{Cond1})).

Further, we multiply (\ref{Equality}) by $\dfrac{du}{dt}$ and use Corollary 1
in \cite{CVL19} to obtain that%
\begin{align}
\left\Vert \frac{du}{dt}\left(  t\right)  \right\Vert ^{2}+\frac{1}{2}\frac
{d}{dt}\left\Vert \dfrac{\partial u}{\partial x}\right\Vert ^{2}  &  \leq
b_{1}\left\Vert \dfrac{du}{dt}\left(  t\right)  \right\Vert +\omega
_{1}\left\Vert u\left(  t\right)  \right\Vert \left\Vert \dfrac{du}{dt}\left(
t\right)  \right\Vert \label{IneqDeriv}\\
&  \leq b_{1}^{2}+\omega_{1}^{2}\left\Vert u\left(  t\right)  \right\Vert
^{2}+\frac{1}{2}\left\Vert \dfrac{du}{dt}\left(  t\right)  \right\Vert
^{2}.\nonumber
\end{align}
For $0\leq t-\alpha\leq r\leq t$ we integrate over the interval $\left(
r,t\right)  $. Hence, by (\ref{IneqL2}) we have%
\[
\left\Vert \dfrac{\partial u}{\partial x}(t)\right\Vert ^{2}\leq\left\Vert
\dfrac{\partial u}{\partial x}(r)\right\Vert ^{2}+2b_{1}^{2}+2\omega_{1}%
^{2}e^{-\delta\left(  t-\alpha\right)  }\left\Vert u\left(  0\right)
\right\Vert ^{2}+\frac{2C_{1}\omega_{1}^{2}}{\delta}.
\]
Integrating now with respect to the variable $r$ over the interval $\left(
t-\alpha,t\right)  $ and using (\ref{IneqL2}) and (\ref{IneqIntH1}) we get%
\begin{align}
\alpha\left\Vert \dfrac{\partial u}{\partial x}(t)\right\Vert ^{2}  &
\leq\frac{\pi^{2}C_{1}}{\delta}+\frac{\pi^{2}}{\delta}\left\Vert
u(t-\alpha)\right\Vert ^{2}+2b_{1}^{2}+2\omega_{1}^{2}e^{-\delta\left(
t-\alpha\right)  }\left\Vert u\left(  0\right)  \right\Vert ^{2}+\frac
{2C_{1}\omega_{1}^{2}}{\delta}\nonumber\\
&  \leq C\left(  1+e^{-\delta\left(  t-\alpha\right)  }\left\Vert u\left(
0\right)  \right\Vert ^{2}\right)  , \label{IneqAbsorb}%
\end{align}
where $C>0$ is a constant.

We take $\alpha=1$ and define the set $K=\{v\in V:\left\Vert v\right\Vert
_{V}^{2}\leq2C\}.$ The compact embedding $V\subset H$ implies that $K$ is
relatively compact in $H$. Also, as $K$ is weakly closed in $V$, it is closed
in $H$. Thus, $K$ is compact in $H$. From (\ref{IneqAbsorb}) we obtain that
for any bounded set $B\subset H$ there exists $T(B)$ (independent of
$\sigma\in\Sigma$) such that $\phi(t,\sigma,B)\subset K$ for all $t\geq T(B)$
and any $\sigma\in\Sigma$. Let $\mathbb{K}=K\times\Sigma$, which is compact in
$\mathcal{X}$. Any bounded set $\mathbb{B}\subset\mathcal{X}$ satisfies that
$\mathbb{B}\subset\mathcal{P}_{H}\mathbb{B\times}\Sigma$, where $\mathcal{P}%
_{H}\mathbb{B}$ is bounded in $H$. Then we have%
\[
\Pi(t,\mathbb{B})\subset\Pi(t,\mathcal{P}_{H}\mathbb{B}\times\Sigma
)\subset\cup_{\sigma\in\Sigma}\left(  \phi\left(  t,\sigma,\mathcal{P}%
_{H}\mathbb{B}\right)  ,\theta_{t}\sigma\right)  \subset\mathbb{K}\text{ for
}t\geq T(\mathbb{B}).
\]

\end{proof}

\bigskip

\begin{lemma}
\label{ConvergConvex}Let $h_{n}\rightarrow h$ weakly in $L^{1}\left(
t_{1},t_{2};H\right)  $ and let for a.a. $\left(  t,x\right)  $ there is
$N(t,x)$ such that $h_{n}(t,x)$ belong to the closed convex set $C\left(
t,x\right)  $ for all $n\geq N.$ Then $h(t,x)\in C\left(  t,x\right)  $ for
a.a. $\left(  t,x\right)  .$
\end{lemma}

\begin{proof}
By \cite[Proposition 1.1]{Tolstonogov} for a.a. $t\in\left(  t_{1}%
,t_{2}\right)  $ there is a sequence of convex combinations of $\{h_{n}\left(
t\right)  \}$ given by%
\[
y_{n}(t)=\sum_{i=1}^{M_{n}}\lambda_{i}h_{k_{i}}(t),\ \sum_{i=1}^{M_{n}}%
\lambda_{i}=1,\ k_{i}\geq n,
\]
such that $y_{n}(t)\rightarrow h\left(  t\right)  $ in $H$. Since for a.a.
$\left(  t,x\right)  $ $h_{k_{i}}(t,x)\in C(t,x)$ if $n\geq N(t,x)$, we get by
the convexity of $C(t,x)$ that $h\left(  t,x\right)  \in C(t,x)$ for a.a.
$\left(  t,x\right)  $.
\end{proof}

\begin{lemma}
\label{ConvergSolutions}Let (\ref{Cond1}), (\ref{CondAttr2}), (\ref{CondAttr3}%
) hold. Let $u_{0}^{n}\rightarrow u_{0}$ in $H$ and $\sigma_{n}=\left(
b_{n},\omega_{n}\right)  \rightarrow\sigma=\left(  \overline{b},\overline
{\omega}\right)  $ in $\Sigma$. Then for any $u_{n}\in\mathcal{D}%
_{0,\sigma_{n}}(u_{0}^{n})$ there exists a subsequence $u_{n^{\prime}}$ and
$u\in\mathcal{D}_{0,\sigma}(u_{0})$ such that%
\[
u_{n^{\prime}}\rightarrow u\text{ in }C([0,+\infty),H).
\]

\end{lemma}

\begin{proof}
By (\ref{IneqAbsorb}) we have that for any $\varepsilon>0$ there exists
$C_{\varepsilon}>0$ such that%
\begin{equation}
\left\Vert u_{n}\left(  t\right)  \right\Vert _{V}\leq C_{\varepsilon}\text{
for all }t\geq\varepsilon\text{ and any }n. \label{IneqV}%
\end{equation}
We fix an arbitrary $T>\varepsilon$. Integrating over $\left(  \varepsilon
,T\right)  $ in (\ref{IneqDeriv}) and using (\ref{IneqV}) we obtain the
existence of $D_{1}=D_{1}(\varepsilon,T)$ such that%
\begin{equation}
\int_{\varepsilon}^{T}\left\Vert \frac{du}{dt}\right\Vert ^{2}dt\leq D_{1}.
\label{IneqDerivU}%
\end{equation}
In view of (\ref{IneqV}) and the compact embedding $V\subset H$, the sequence
$u_{n}\left(  t\right)  $ is relatively compact in $H$ for all $t\in
$[$\varepsilon,T]$. Also, it follows easily from (\ref{IneqDerivU}) that the
functions $u_{n}:[\varepsilon,T]\rightarrow H$ are equicontinuous. Hence,
(\ref{IneqV}), (\ref{IneqDerivU}) and the Ascoli-Arzel\`{a} theorem imply that
for any $0<\varepsilon<T$ up to a subsequence the following convergences hold:%
\[
u_{n}\rightarrow u\text{ weakly star in }L^{\infty}(\varepsilon,T;V),
\]%
\[
\frac{du_{n}}{dt}\rightarrow\frac{du}{dt}\text{ weakly in }L^{2}%
(\varepsilon,T;H),
\]%
\[
u_{n}\rightarrow u\text{ in }C([\varepsilon,T],H).
\]
Since $u_{n}\left(  \text{\textperiodcentered}\right)  $ satisfies
(\ref{Equality}), where $r_{n}\in L_{loc}^{2}(0,+\infty;H)$ is such that
$r_{n}(t,x)\in b_{n}(t)H_{0}(u_{n}(t,x))+\omega_{n}(t)u_{n}(t,x)$ for a.a.
$\left(  t,x\right)  $, there exists a sequence $h_{n}\in L_{loc}%
^{2}(0,+\infty;H)$ such that $h_{n}(t,x)\in H_{0}(u_{n}(t,x)),$ for a.a.
$\left(  t,x\right)  ,$ and%
\[
-Au_{n}(t)=\frac{du_{n}}{dt}(t)-b_{n}(t)h_{n}(t)-\omega_{n}(t)u_{n}(t)\text{
for a.a. }t.
\]
As up to a subsequence $h_{n}$ converges weakly in $L_{loc}^{2}(0,+\infty;H)$
to some $h$, we have%
\[
-Au_{n}\rightarrow\frac{du}{dt}-\overline{b}h-\overline{\omega}u\text{ weakly
in }L^{2}(\varepsilon,T;H)\text{ for all }0<\varepsilon<T.
\]
Hence, as $Au_{n}=-\dfrac{\partial^{2}u_{n}}{\partial x^{2}}$ converges to
$Au=-\dfrac{\partial^{2}u}{\partial x^{2}}$ in the sense of distributions, we
have%
\[
\frac{du}{dt}+Au\left(  t\right)  =\overline{b}\left(  t\right)  h\left(
t\right)  +\overline{\omega}\left(  t\right)  u\left(  t\right)  \text{ for
a.a. }t\in\left(  0,+\infty\right)  .
\]

We need to show that $h\left(  t,x\right)  \in H_{0}(u(t,x))$ for a.a.
$\left(  t,x\right)  $. We observe that for a.a. $\left(  t,x\right)  $ there
is $N(t,x)$ such that $h_{n}(t,x)\in H_{0}(u(t,x)))$ if $n\geq N$. Indeed, let
$\left(  t_{0},x_{0}\right)  \in A^{c}$, where $A$ is a set of measure $0$
such that $u_{n}\left(  t,x\right)  \rightarrow u\left(  t,x\right)  $ for any
$\left(  t,x\right)  \in A^{c}$. If $u(t_{0},x_{0})>0$ ($<0$), then
$u_{n}(t_{0},x_{0})>0$ ($<0$) for $n$ large enough, so $h_{n}(t_{0},x_{0})=1$
($-1$). Hence, $h_{n}(t_{0},x_{0})\in H_{0}(u(t,x))$. If $u\left(  t_{0}%
,x_{0}\right)  =0$, then $h_{n}(t_{0},x_{0})\in\lbrack-1,1]=H_{0}(u(t,x)).$
Then the assertion follows from Lemma \ref{ConvergConvex}.

It remains to show that $u\left(  t\right)  \rightarrow u\left(  0\right)  $
as $t\rightarrow0^{+}$ and that $u_{n}(t_{n})\rightarrow u\left(  0\right)  $
as $t_{n}\rightarrow0$.

Let $v_{n}(t)=u_{n}(t)-\widehat{u}\left(  t\right)  $, where $\widehat{u}$ is
the unique solution of the linear problem%
\[
\left\{
\begin{array}
[c]{c}%
\dfrac{\partial u}{\partial t}-\dfrac{\partial^{2}u}{\partial x^{2}}%
=\overline{\omega}\left(  t\right)  u,\\
u\left(  t,0\right)  =u(t,1)=0,\\
u\left(  0,x\right)  =u_{0}(x).
\end{array}
\right.
\]
Fix some $T>0$. Then in a standard way we obtain that%
\begin{align*}
\frac{1}{2}\frac{d}{dt}\left\Vert v_{n}\right\Vert ^{2}+\left\Vert
v_{n}\left(  t\right)  \right\Vert _{V}^{2}  &  \leq b_{n}\left(  t\right)
\left\Vert v_{n}\left(  t\right)  \right\Vert +\omega_{n}\left(  t\right)
\left\Vert v_{n}\left(  t\right)  \right\Vert ^{2}+\left\vert \omega
_{n}(t)-\overline{\omega}\left(  t\right)  \right\vert \left\Vert
\widehat{u}\left(  t\right)  \right\Vert \left\Vert v_{n}\left(  t\right)
\right\Vert \\
&  \leq\frac{b_{1}^{2}}{2\varepsilon_{0}}+\left(  \omega_{1}+\varepsilon
_{0}\right)  \left\Vert v_{n}\left(  t\right)  \right\Vert ^{2}+\alpha
_{n},\text{ for a.a. \ }t\in\left(  0,T\right)  ,
\end{align*}
where $\varepsilon_{0}>0$ is such that $\omega_{1}+\varepsilon_{0}<\pi^{2}$
and $\alpha_{n}=\frac{1}{2\varepsilon_{0}}\sup_{t\in\lbrack0,T]}\left(
\left\vert \omega_{n}(t)-\overline{\omega}\left(  t\right)  \right\vert
^{2}\left\Vert \widehat{u}\left(  t\right)  \right\Vert ^{2}\right)
\rightarrow0$ as $n\rightarrow\infty$. Then%
\[
\left\Vert v_{n}\left(  t\right)  \right\Vert ^{2}\leq\left\Vert v_{n}\left(
0\right)  \right\Vert ^{2}+Ct,
\]
where $\frac{b_{1}^{2}}{\varepsilon_{0}}+2\alpha_{n}\leq C$ for any $n$. It
follows that $\left\Vert u\left(  t\right)  -\widehat{u}(t)\right\Vert
^{2}=\lim_{n\rightarrow\infty}\left\Vert v_{n}\left(  t\right)  \right\Vert
^{2}\leq Ct$ for $t>0.$ Thus,%
\[
\left\Vert u\left(  t\right)  -u_{0}\right\Vert ^{2}\leq2\left\Vert u\left(
t\right)  -\widehat{u}(t)\right\Vert ^{2}+2\left\Vert \widehat{u}%
(t)-u_{0}\right\Vert ^{2}\rightarrow0\text{ as }t\rightarrow0^{+}.
\]
Therefore, $u\in\mathcal{D}_{0,\sigma}(u_{0})$.

Finally, if $t_{n}\rightarrow0^{+}$, then%
\begin{align*}
\left\Vert u_{n}\left(  t_{n}\right)  -u_{0}\right\Vert ^{2}  &
\leq2\left\Vert v_{n}\left(  t_{n}\right)  \right\Vert ^{2}+2\left\Vert
\widehat{u}(t_{n})-u_{0}\right\Vert ^{2}\\
&  \leq2\left\Vert v_{n}\left(  0\right)  \right\Vert ^{2}+2Ct_{n}+2\left\Vert
\widehat{u}(t_{n})-u_{0}\right\Vert ^{2}\rightarrow0,
\end{align*}
so $u_{n}\rightarrow u$ in $C([0,+\infty),H).$
\end{proof}

\bigskip

For any $\sigma\in\Sigma$ and $u\in\mathcal{R}_{0,\sigma}$ let $y:\mathbb{R}%
^{+}\rightarrow\mathcal{X}$ be given by $y\left(  \text{\textperiodcentered
}\right)  =\left(  u\left(  \text{\textperiodcentered}\right)  ,\theta
_{\text{\textperiodcentered}}\sigma\right)  $ and denote by $\mathcal{K}%
\subset C([0,+\infty),\mathcal{X})$ the set of all functions $y\left(
\text{\textperiodcentered}\right)  $ of this type. We consider the following
standard axiomatic properties:

\begin{itemize}
\item[$\left(  K1\right)  $] For any $y_{0}=\left(  u_{0},\sigma\right)
\in\mathcal{X}$ there exists $y\in\mathcal{K}$ such that $y\left(  0\right)
=y_{0}.$

\item[$\left(  K2\right)  $] For any $y\in\mathcal{K}$ and $s\geq0$,
$y_{s}=y($\textperiodcentered$+s)\in\mathcal{K}.$

\item[$\left(  K3\right)  $] If $y_{1},y_{2}\in\mathcal{K}$ are such that
$y_{2}(0)=y_{1}(s)$, then the composition%
\[
y\left(  t\right)  =\left\{
\begin{array}
[c]{c}%
y_{1}(t)\text{ if }0\leq t\leq s,\\
y_{2}(t-s)\text{ if }t\geq s,
\end{array}
\right.
\]
belongs to $\mathcal{K}$.

\item[$\left(  K4\right)  $] If $y_{n}\in\mathcal{K}$ is such that
$y_{n}(0)\rightarrow y_{0}$, then there exists a subsequence $\{y_{n^{\prime}%
}\}$ and $y\in\mathcal{K}$ such that $y_{n^{\prime}}(t)\rightarrow y(t)$ in
$\mathcal{X}$ uniformly in compact sets of $[0,+\infty)$.
\end{itemize}

We recall that the map $\Pi\left(  t,\text{\textperiodcentered}\right)
:\mathcal{X\rightarrow X}$ is upper semicontinuous if for any $u_{0}%
\in\mathcal{X}$ and any neighborhood $O$ of $\Pi(t,u_{0})$ there exists
$\rho>0$ such that $\Pi\left(  t,u\right)  \subset O$ as soon as
$\rho_{\mathcal{X}}\left(  u,u_{0}\right)  <\rho$.

\begin{lemma}
\label{PropertiesPi}Let (\ref{Cond1}), (\ref{CondAttr2}), (\ref{CondAttr3})
hold. Then properties $\left(  K1\right)  -\left(  K4\right)  $ hold. An
alternative definition for the map $\Pi$ is
\begin{equation}
\Pi(t,u_{0})=\{\xi\in\mathcal{X}:\xi=y\left(  t\right)  ,\ y\in\mathcal{K}\}
\label{MapPi}%
\end{equation}
and for any $t\geq0$ the map $\Pi\left(  t,\text{\textperiodcentered}\right)
$ is upper semicontinuous and has compact values.
\end{lemma}

\begin{proof}
$\left(  K1\right)  $ follows from the existence of solutions to (\ref{Incl})
for any $u_{0}\in H$ and $\sigma\in\Sigma$. $\left(  K2\right)  $ follows from%
\[
y(\text{\textperiodcentered}+s)=\left(  u\left(  \text{\textperiodcentered
}+s\right)  ,\theta_{\text{\textperiodcentered}+s}\sigma\right)  =\left(
u_{s}\left(  \text{\textperiodcentered}\right)  ,\theta
_{\text{\textperiodcentered}}\theta_{s}\sigma\right)  \in\mathcal{K}\text{,}%
\]
as $u_{s}\left(  \text{\textperiodcentered}\right)  \in\mathcal{R}%
_{0,\theta_{s}\sigma}$ (see the proof of Lemma \ref{TranslationLemma}). For
$\left(  K3\right)  $ we note that $y_{1}\left(  t\right)  =\left(
u_{1}\left(  t\right)  ,\theta_{t}\sigma_{1}\right)  $, $y_{2}\left(
t-s\right)  =\left(  u_{2}\left(  t-s\right)  ,\theta_{t-s}\sigma_{2}\right)
$, $u_{2}\left(  0\right)  =u_{1}\left(  s\right)  $, $\sigma_{2}=\theta
_{s}\sigma_{1}$, so that $y\left(  t\right)  =\left(  u\left(  t\right)
,\theta_{t}\sigma_{1}\right)  $, for any $t\geq0$, where
\[
u\left(  t\right)  =\left\{
\begin{array}
[c]{c}%
u_{1}(t)\text{ if }0\leq t\leq s,\\
u_{2}(t-s)\text{ if }t\geq s.
\end{array}
\right.
\]
If we prove that $u\in\mathcal{R}_{0,\sigma_{1}}$, then $y\in\mathcal{K}$.
Indeed, if $\sigma_{1}\left(  t\right)  =\left(  b_{1}\left(  t\right)
,\omega_{1}\left(  t\right)  \right)  $, then for $t\geq s$ we have%
\begin{align*}
u\left(  t\right)   &  =u_{2}\left(  t-s\right)  =e^{-A(t-s)}u_{1}\left(
s\right)  +\int_{0}^{t-s}e^{-A(t-s-\tau)}(b_{1}\left(  \tau+s\right)
h_{2}(\tau)+\omega_{1}\left(  \tau+s\right)  u_{2}(\tau))d\tau\\
&  =e^{-A(t-s)}e^{-As}u_{1}\left(  0\right)  +e^{-A(t-s)}\int_{0}%
^{s}e^{-A(s-\tau)}(b_{1}\left(  \tau\right)  h_{1}(\tau)+\omega_{1}\left(
\tau\right)  u_{1}(\tau))d\tau\\
&  +\int_{s}^{t}e^{-A(t-\tau)}(b_{1}\left(  \tau\right)  h_{2}(\tau
-s)+\omega_{1}\left(  \tau\right)  u_{2}(\tau-s))d\tau\\
&  =e^{-At}u_{1}\left(  0\right)  +\int_{0}^{t}e^{-A\left(  t-\tau\right)
}\left(  (b_{1}\left(  \tau\right)  h(\tau)+\omega_{1}\left(  \tau\right)
u(\tau))d\tau\right)  ,
\end{align*}
where $h_{i}\in L_{loc}^{2}(0,+\infty;H)$ are such that $h_{i}\left(
\tau\right)  \in H_{0}(u_{i}(\tau))$ for a.a. $\tau>0$, and
\[
h\left(  \tau\right)  =\left\{
\begin{array}
[c]{c}%
h_{1}(t)\text{ if }0\leq\tau\leq s,\\
h_{2}(t-s)\text{ if }\tau\geq s,
\end{array}
\right.
\]
belongs to $L_{loc}^{2}(0,+\infty;H)$ and satisfies $h\left(  \tau\right)  \in
H_{0}(u(\tau))$ for a.a. $\tau>0$. Hence, $u$ is a mild solution for the
symbol $\sigma_{1}$.

Property $\left(  K4\right)  $ follows from Lemma \ref{ConvergSolutions} and
implies easily that $\Pi\left(  t,\text{\textperiodcentered}\right)  $ has
compact values and is upper semicontinuous.
\end{proof}

\bigskip

We are now in position of establishing the existence of the global attractor
and its relationship with the cocycle attractor.

\begin{theorem}
\label{AttractorPi}The multivalued semiflow $\Pi$ possesses the global compact
invariant attractor $\mathbb{A}$. Moreover,
\[
\mathbb{A}=\cup_{\sigma\in\Sigma}A(\sigma)\times\{\sigma\},
\]
where $\{\mathcal{A}(\sigma)\}_{\sigma\in\Sigma}$ is the cocycle attractor.
\end{theorem}

\begin{proof}
The existence of the global compact invariant attractor follows from Lemmas
\ref{CompactAbsorbing}, \ref{PropertiesPi} and \cite[Theorem 4 and Remark
8]{MelnikValero98}. The relationship with the cocycle attractor is a
consequence of \cite[Corollary 1]{Cui}.
\end{proof}

\bigskip

As before, we will characterize the global attractor in terms of bounded
complete trajectories. A complete trajectory of $\mathcal{K}$ is a function
$\Phi:\mathbb{R}\rightarrow\mathcal{X}$ such that $\Phi\left(
\text{\textperiodcentered}+s\right)  \mid_{t\geq0}\in\mathcal{K}$ for any
$s\in\mathbb{R}$.

\begin{lemma}
\label{EquivalenceTrajectories}Assume that (\ref{Cond1}), (\ref{CondAttr2}),
(\ref{CondAttr3}) hold. Let $\gamma$ be a complete trajectory of
$\mathcal{R}_{\sigma}$. Then%
\begin{equation}
\Phi\left(  t\right)  =\left(  \gamma\left(  t\right)  ,\theta_{t}%
\sigma\right)  \text{, for any }t\in\mathbb{R}\text{,} \label{EqCompleteTraj}%
\end{equation}
is a complete trajectory of $\mathcal{K}$. Conversely, if $\Phi$ is a complete
trajectory of $\mathcal{K}$, then there exist $\sigma\in\Sigma$ and a complete
trajectory $\gamma$ of $\mathcal{R}_{\sigma}$ such that (\ref{EqCompleteTraj}) holds.
\end{lemma}

\begin{proof}
Let $\gamma$ be a complete trajectory of $\mathcal{R}_{\sigma}$. Since
$\gamma\left(  \text{\textperiodcentered}+s\right)  \in\mathcal{D}%
_{0,\theta_{s}\sigma}(\gamma\left(  s\right)  )$, for any $s\in\mathbb{R}$, we
have that $y\left(  \text{\textperiodcentered}\right)  =(\gamma($%
\textperiodcentered$+s),\theta_{\text{\textperiodcentered}}\theta_{s}%
\sigma)\in\mathcal{K}$. Thus, $\Phi$ is a complete trajectory of $\mathcal{K}$.

Further, let $\Phi=\left(  \gamma,\alpha\right)  $ be a complete trajectory of
$\mathcal{K}$. Then%
\[
\left(  \gamma\left(  \text{\textperiodcentered}+s\right)  ,\alpha\left(
\text{\textperiodcentered}+s\right)  \right)  \mid_{t\geq0}\in\mathcal{K}%
\text{ for any }s\text{,}%
\]
so $\gamma\left(  \text{\textperiodcentered}+s\right)  \in\mathcal{R}%
_{0,\alpha\left(  s\right)  }=\mathcal{R}_{0,\theta_{s}\alpha\left(  0\right)
}$. Therefore, as explained in the previous subsection, $\gamma$ is a complete
trajectory of $\mathcal{R}_{\sigma}$ with $\sigma=\alpha\left(  0\right)  $.
\end{proof}

\bigskip

\begin{theorem}
\label{AttrCharac}Let (\ref{Cond1}), (\ref{CondAttr2}), (\ref{CondAttr3})
hold. Then the global attractor $\mathbb{A}$ is given by%
\begin{align}
\mathbb{A}  &  =\{\Phi\left(  0\right)  :\Phi\text{ is a bounded complete
trajectory of }\mathcal{K}\}\nonumber\\
&  =\cup_{t\in\mathbb{R}}\{\Phi\left(  t\right)  :\Phi\text{ is a bounded
complete trajectory of }\mathcal{K}\}\nonumber\\
&  =\cup_{\sigma\in\Sigma}\{\left(  \gamma\left(  0\right)  ,\sigma\right)
:\gamma\text{ is a bounded complete trajectory of }\mathcal{R}_{\sigma
}\}\nonumber\\
&  =\cup_{\sigma\in\Sigma,t\in\mathbb{R}}\{\left(  \gamma\left(  t\right)
,\theta_{t}\sigma\right)  :\gamma\text{ is a bounded complete trajectory of
}\mathcal{R}_{\sigma}\}. \label{CharacterizationA}%
\end{align}

\end{theorem}

\begin{proof}
We obtain the first two equalities from \cite[Theorems 9 or 10]{KKV14},
whereas the other ones follow from Lemma \ref{EquivalenceTrajectories}.
\end{proof}

\bigskip

When we restrict the solutions to the positive cone, we obtain a Morse
structure of the global attractor.

We denote by $\mathcal{X}^{+}$ the product space $H^{+}\times\Sigma$. Let now
take the skew product semiflow $\Pi^{+}:\mathbb{R}^{+}\times\mathcal{X}%
^{+}\rightarrow\mathcal{P}(\mathcal{X}^{+})$ given by
\begin{align*}
\Pi^{+}(t,\left(  x,\sigma\right)  )  &  =\left(  \phi^{+}\left(
t,\sigma,x\right)  ,\theta_{t}\sigma\right) \\
&  =\{\xi\in\mathcal{X}^{+}:\xi=y\left(  t\right)  ,\ y\in\mathcal{K}^{+}\},
\end{align*}
where $\mathcal{K}^{+}$ stands for the subset of $C([0,+\infty),\mathcal{X}%
^{+})$ of functions $y\left(  \text{\textperiodcentered}\right)  =\left(
u\left(  \text{\textperiodcentered}\right)  ,\theta_{\text{\textperiodcentered
}}\sigma\right)  $ such that $\sigma\in\Sigma$ and $u\in\mathcal{R}_{0,\sigma
}^{+}$. The strict multivalued semiflow $\Pi^{+}$ possesses the global compact
invariant attractor $\mathbb{A}^{+}=\mathbb{A}\cap\mathcal{X}^{+}$, which
satisfies%
\begin{align}
\mathbb{A}^{+}  &  =\{\Phi\left(  0\right)  :\Phi\text{ is a bounded complete
trajectory of }\mathcal{K}^{+}\}\nonumber\\
&  =\cup_{t\in\mathbb{R}}\{\Phi\left(  t\right)  :\Phi\text{ is a bounded
complete trajectory of }\mathcal{K}^{+}\}\nonumber\\
&  =\cup_{\sigma\in\Sigma}\{\left(  \gamma\left(  0\right)  ,\sigma\right)
:\gamma\text{ is a bounded complete trajectory of }\mathcal{R}_{\sigma}%
^{+}\}\nonumber\\
&  =\cup_{\sigma\in\Sigma,t\in\mathbb{R}}\{\left(  \gamma\left(  t\right)
,\theta_{t}\sigma\right)  :\gamma\text{ is a bounded complete trajectory of
}\mathcal{R}_{\sigma}^{+}\}. \label{CharacterizationA+}%
\end{align}

A set $M$ is called weakly invariant if for any $y_{0}\in M$ there exist a
complete trajectory $\Phi$ of $\mathcal{K}^{+}$ such that $\cup_{t\in
\mathbb{R}}\Phi\left(  t\right)  \subset M$ and $\phi\left(  0\right)  =y_{0}%
$. It is obvious that the global attractor $\mathbb{A}^{+}$ is weakly
invariant. A weakly invariant set $M$ is said to be isolated if there exists a
neighborhood $O$ of $M$ such that $M$ is the maximal weakly invariant set in
it. If $M$ is compact, this is equivalent to saying that there exists an
$\varepsilon>0$ such that $M$ is the maximal weakly invariant set in the
$\varepsilon$-neighborhood $O_{\varepsilon}(M)=\{y\in\mathcal{X}%
^{+}:dist_{\mathcal{X}^{+}}(y,M)<\varepsilon\}$. As we will consider weakly
invariant sets $M$ belonging to the global attractor $\mathbb{A}^{+}$, they
will be necessarily compact \cite[Lemma 19]{CostaValero16}.

A family of isolated weakly invariant sets $\mathcal{M}=\{M_{1},...,M_{n}\}$
in $\mathbb{A}^{+}$ is called disjoint if there is $\delta>0$ such that
$O_{\delta}(M_{i})\cap O_{\delta}(M_{j})=\varnothing$ for all $i\not =j$.

We say that the multivalued semiflow $\Pi^{+}$ is dynamically gradient with
respect to the disjoint family of isolated weakly invariant sets
$\mathcal{M}=\{M_{1},...,M_{n}\}$ in $\mathbb{A}^{+}$ if for any bounded
complete trajectory $\Phi$ of $\mathcal{K}^{+}$ we have that either
$\cup_{t\in\mathbb{R}}\Phi\left(  t\right)  \in M_{i}$ for some $i\in
\{1,...,n\}$ or
\begin{align}
dist_{\mathcal{X}^{+}}(\Phi\left(  t\right)  ,M_{i})  &  \rightarrow0\text{ as
}t\rightarrow+\infty,\label{Converg+}\\
dist_{\mathcal{X}^{+}}(\Phi\left(  t\right)  ,M_{j})  &  \rightarrow0\text{ as
}t\rightarrow-\infty, \label{Converg-}%
\end{align}
for some $1\leq i<j\leq n$. We observe that convergences (\ref{Converg+}) and
(\ref{Converg-}) are equivalent to saying that $\omega\left(  \Phi\right)
\subset M_{i}$ and $\alpha\left(  \Phi\right)  \subset M_{j}$, where
$\omega\left(  \Phi\right)  $, $\alpha\left(  \Phi\right)  $ stand,
respectively, for the omega and alpha limit sets of the bounded complete
trajectory $\Phi$.

For each $\sigma$ we denote by $\xi_{M,\sigma}$ the unique positive bounded
complete trajectory (the positive nonautonomous equilibrium) given in Section
\ref{Nonautonomous}. We define then the following compact weakly invariant
sets in $\mathbb{A}^{+}$:%
\begin{align*}
M_{1}  &  =\{\left(  \xi_{M,\sigma}\left(  0\right)  ,\sigma\right)
:\sigma\in\Sigma\}\\
&  =\{\left(  \xi_{M,\sigma}\left(  t\right)  ,\theta_{t}\sigma\right)
:\sigma\in\Sigma,\ t\in\mathbb{R}\},\\
M_{2}  &  =\{0\}.
\end{align*}

Since $\xi_{M,\sigma}\left(  0\right)  \geq w_{b_{0},\omega_{0}}^{+}$ for any
$\sigma$, these sets are clearly disjoint. From the results in Section
\ref{Nonautonomous} and (\ref{CharacterizationA+}) we can see that, apart from
$M_{1}$ and $M_{2}$, the only elements in the global attractor $\mathbb{A}%
^{+}$ are $\left(  \gamma\left(  t\right)  ,\theta_{t}\sigma\right)  $, where
$\gamma$ is a bounded complete trajectory of $\mathcal{R}_{\sigma}^{+}$ of the
type (\ref{Connection3}) with $u\left(  \text{\textperiodcentered}\right)  $
being the solution to the linear problem (\ref{LinearNonAut}) for the symbol
$\sigma=\left(  \overline{b},\overline{\omega}\right)  $. We know that
\[
\left\Vert \gamma\left(  t\right)  -\xi_{M,\sigma}\left(  t\right)
\right\Vert \rightarrow0\text{ as }t\rightarrow+\infty,
\]
so%
\begin{align*}
dist\left(  \left(  \gamma\left(  t\right)  ,\theta_{t}\sigma\right)
,M_{1}\right)   &  \rightarrow0\text{ as }t\rightarrow+\infty,\\
dist\left(  \left(  \gamma\left(  t\right)  ,\theta_{t}\sigma\right)
,M_{2}\right)   &  \rightarrow0\text{ as }t\rightarrow-\infty.
\end{align*}
The sets $M_{1},M_{2}$ are isolated. Indeed, we take disjoint $\varepsilon
$-neighborhoods of these sets, $O_{\varepsilon}\left(  M_{1}\right)  $ and
$O_{\varepsilon}\left(  M_{2}\right)  $. In $O_{\varepsilon}\left(
M_{1}\right)  $ the only possible complete trajectories are $\xi_{M,\sigma}$,
because any other one should converge to $0$ as $t\rightarrow-\infty$ and then
leave the neighborhood $O_{\varepsilon}\left(  M_{1}\right)  $. By the same
reason the only possible bounded complete trajectory in $O_{\varepsilon
}\left(  M_{2}\right)  $ is $0$, so that $M_{1}$ and $M_{2}$ are the maximal
weakly invariant sets in $O_{\varepsilon}\left(  M_{1}\right)  $ and
$O_{\varepsilon}\left(  M_{2}\right)  $, respectively.

All in all, we have shown the following.

\begin{theorem}
Let (\ref{Cond1}), (\ref{CondAttr2}), (\ref{CondAttr3}) hold. Then $\Pi^{+}$
is dynamically gradient with respect to the disjoint family of isolated weakly
invariant sets $\mathcal{M}=\{M_{1},M_{2}\}.$ Hence, the global attractor
$\mathbb{A}^{+}$ possesses a gradient structure.
\end{theorem}

\subsection{The uniform attractor}

We finish this section by showing the relationship of the previous attractors
with the uniform attractor of the cocycle $\phi$.

It follows from Theorem 5 in \cite{Cui} that $\mathcal{A}=\mathcal{P}%
_{H}\mathbb{A}$ ($\mathbb{A}$ is the attractor of the skew product flow) is
the uniform attractor for $\phi$, which means that:

\begin{itemize}
\item $\mathcal{A}$ is compact.

\item $\mathcal{A}$ is uniformly attracting, that is, for any bounded set
$B\subset H$ we have%
\begin{equation}
\sup_{\sigma\in\Sigma}dist\left(  \phi\left(  t,\sigma,B\right)
,\mathcal{A}\right)  \rightarrow0\text{, as }t\rightarrow+\infty.
\label{AttractionUnif}%
\end{equation}

\item $\mathcal{A}$ is the minimal closed set satisfying (\ref{AttractionUnif}).
\end{itemize}

By characterization (\ref{CharacterizationA}) we obtain that
\begin{align*}
\mathcal{A}  &  \mathcal{=}\cup_{\sigma\in\Sigma}\{\gamma\left(  0\right)
:\gamma\text{ is a bounded complete trajectory of }\mathcal{R}_{\sigma}\}\\
&  =\cup_{\sigma\in\Sigma,t\in\mathbb{R}}\{\gamma\left(  t\right)
:\gamma\text{ is a bounded complete trajectory of }\mathcal{R}_{\sigma}\}.
\end{align*}
This implies, in particular, that $\mathcal{A}$ is weakly invariant.

Also, by Theorem 8 in \cite{Cui} we have the relationship of the uniform
attractor with the cocycle and pullback attractors:%
\[
\mathcal{A}=\cup_{\sigma\in\Sigma}\mathcal{A}(\sigma)=\cup_{\sigma\in\Sigma
}\mathcal{A}_{\sigma}(0).
\]

In the same way, $\mathcal{A}^{+}=\mathcal{P}_{H}\mathbb{A}^{+}$ is the
uniform attractor for $\phi^{+}$ and
\begin{align*}
\mathcal{A}^{+}  &  \mathcal{=}\cup_{\sigma\in\Sigma}\{\gamma\left(  0\right)
:\gamma\text{ is a bounded complete trajectory of }\mathcal{R}_{\sigma}%
^{+}\}\\
&  =\cup_{\sigma\in\Sigma,t\in\mathbb{R}}\{\gamma\left(  t\right)
:\gamma\text{ is a bounded complete trajectory of }\mathcal{R}_{\sigma}^{+}\},
\end{align*}%
\[
\mathcal{A}^{+}=\cup_{\sigma\in\Sigma}\mathcal{A}^{+}(\sigma)=\cup_{\sigma
\in\Sigma}\mathcal{A}_{\sigma}^{+}(0).
\]
Let $\mathcal{C}=\cup_{\sigma\in\Sigma}\cup_{\gamma_{\sigma}\in\mathcal{L}%
_{\sigma}}\gamma_{\sigma}\left(  0\right)  $, where $\mathcal{L}_{\sigma}$ is
the set of all bounded complete trajectory of $\mathcal{R}_{\sigma}^{+}$ of
the type (\ref{Connection3}) with $u\left(  \text{\textperiodcentered}\right)
$ being the solution to the linear problem (\ref{LinearNonAut}) for the symbol
$\sigma=\left(  \overline{b},\overline{\omega}\right)  $. Hence,%
\[
\mathcal{A}^{+}=\{0\}\cup\{\xi_{M,\sigma}\left(  0\right)  :\sigma\in
\Sigma\}\cup\mathcal{C}.
\]

\section{Appendix}

In this appendix we present some auxiliary results that are necessary for the
arguments throughout this paper.

\subsection{A maximum principle for non-smooth functions}

Usually in the literature the maximum principle is stated for smooth functions
(see for example \cite{Protter}). However, we need in this paper a maximum
principle for less regular functions. Such result is proved in \cite{Kadlec}
in a rather general setting. We describe in this appendix a particular
situation which is derived from the theorems in \cite{Kadlec}.

Let $\mathcal{O}$ be a region in $\mathbb{R}^{2}$ and let $\left(  t_{0}%
,x_{0}\right)  \in\mathcal{O}$ and $\rho,\sigma>0.$ We denote%
\[
Q_{\rho,\sigma}=\{(t,x):t\in(t_{0}-\sigma,t_{0}),\left\vert x-x_{0}\right\vert
<\rho\},
\]
where we assume that $t_{0},x_{0},\rho,\sigma$ are such that $\overline
{Q}_{\rho,\sigma}\subset\mathcal{O}$.

We denote by $W$ the space of all functions from $L^{2}\left(  \mathcal{O}%
\right)  $ such that%
\[
\int_{\mathcal{O}}\left(  \left\vert u\left(  t,x\right)  \right\vert
^{2}+\left\vert \frac{\partial u}{\partial x}\left(  t,x\right)  \right\vert
^{2}\right)  d\mu<+\infty.
\]

As a particular case of Theorem 6.4 in \cite{Kadlec} we obtain the following
maximum and minimum principles.

\begin{theorem}
\label{Maximum}(Maximum principle) Let $u\in W$ be such that
\[
\frac{\partial u}{\partial t}-\frac{\partial^{2}u}{\partial x^{2}}\leq0
\]
in the sense of distributions. If%
\[
ess\sup{}_{\left(  t,x\right)  \in Q_{\rho\nu,\sigma_{1}}}u(t,x)=M,
\]
for some $\nu$, $0<\nu<1$, and any $\sigma_{1}$, where $0<\sigma_{1}<\sigma$,
then $u\left(  t,x\right)  =M$ for a.a. $\left(  t,x\right)  \in
Q_{\rho,\sigma}.$
\end{theorem}

\begin{theorem}
\label{Minimum}(Minimum principle) Let $u\in W$ be such that
\[
\frac{\partial u}{\partial t}-\frac{\partial^{2}u}{\partial x^{2}}\geq0
\]
in the sense of distributions. If%
\[
ess\inf{}_{\left(  t,x\right)  \in Q_{\rho\nu,\sigma_{1}}}u(t,x)=M,
\]
for some $\nu$, $0<\nu<1$, and any $\sigma_{1}$, where $0<\sigma_{1}<\sigma$,
then $u\left(  t,x\right)  =M$ for a.a. $\left(  t,x\right)  \in
Q_{\rho,\sigma}.$
\end{theorem}

\subsection{Parabolic regularity}

The regularity of solutions of the linear parabolic problem is well known in
the literature (see e.g. \cite{AgranovichVishik}, \cite{Bardos},
\cite{Lions61}, \cite{SellBook}). In this section, we consider the linear
parabolic problem (\ref{Aux}) and following the result given in Theorem 42.14
in \cite{SellBook} we describe some regularity results of its solutions in a
suitable form for our purposes in this paper.

The following proposition is proved in the same way as Theorem 42.14 in
\cite{SellBook}. However, as the assumptions on the function $f$ do not
coincide exactly, we give a sketch of the proof.

\begin{proposition}
\label{Regularity}Let $f\in L_{loc}^{2}\left(  \tau,+\infty;H\right)  \cap
W_{loc}^{1,2}(\tau+\varepsilon,+\infty;H)$ for all $\varepsilon>0$. For any
$u_{\tau}\in H$, the solution $u\left(  \text{\textperiodcentered}\right)  $
to problem (\ref{Aux}) belongs to $C^{1}((\tau,+\infty),V)$ and $C^{0,\frac
{1}{2}}([\tau+\varepsilon,T],H^{2}(0,1)),\ \dfrac{d^{2}u}{dt^{2}}\in
L^{2}(\tau+\varepsilon,T;H)$ for all $\tau<T$ and $\varepsilon>0$.

If, moreover, $f\in C((\tau,+\infty),H^{1}(0,1))$, then $u\in C((\tau
,+\infty),H^{3}(0,1)$.
\end{proposition}

\begin{proof}
We fix an arbitrary $T>\tau$. We know by Lemma \ref{RegularityA} and Remark
\ref{RegularityB} that $u\in C((\tau,T],V)$ and $u\in L^{2}(\tau
+\varepsilon,T;D(A))$ for any $\varepsilon>0.$

We divide the proof into three steps.

\textbf{Step 1. }Let $u_{\tau}\in D(A)$ and $f\in W_{loc}^{1,2}(\tau
,+\infty;H).$

We put $v_{\tau}=f\left(  \tau\right)  -Au_{\tau}\in H$ and $g=\dfrac{df}%
{dt}\in L^{2}(\tau,T;H)$. Denote by $\omega\left(  \text{\textperiodcentered
}\right)  $ the unique strong solution of problem
\begin{equation}
\left\{
\begin{array}
[c]{l}%
\dfrac{\partial\omega}{\partial t}-\dfrac{\partial^{2}\omega}{\partial x^{2}%
}=g(t),\text{ on }(\tau,+\infty)\times(0,1),\\
\omega(t,0)=\omega(t,1)=0,\\
\omega(\tau,x)=v_{\tau}(x),
\end{array}
\right.  \label{Aux2}%
\end{equation}
Again, by Lemma \ref{RegularityA} and Remark \ref{RegularityB} we have that
$\omega\in C((\tau,T],V)$, $\omega\in L^{2}(\tau+\varepsilon,T;D(A)),\ \dfrac
{d\omega}{dt}\in L^{2}(\tau+\varepsilon,T;H)$ for any $\varepsilon>0$.

Let $z(t)=u_{\tau}+\int_{\tau}^{t}\omega(s)ds$. Hence, $z\in C^{1}%
((\tau,T],V)$ and%
\[
\left\Vert z(t)-z(s)\right\Vert _{H^{2}}\leq\int_{s}^{t}\left\Vert
\omega(r)\right\Vert _{H^{2}}dr\leq(t-s)^{\frac{1}{2}}\int_{s}^{t}\left\Vert
\omega(r)\right\Vert _{H^{2}}^{2}dr\text{ for any }\tau<s<t,
\]
so $z\in C^{0,\frac{1}{2}}([\tau+\varepsilon,T],H^{2}(0,1))$. Following the
proof of Theorem 42.14 in \cite{SellBook} we obtain that $z$ is a strong
solution to problem (\ref{Aux}) on $\left(  \tau,T\right)  $. By uniqueness of
solutions $u=z$, so
\begin{equation}
u\in C^{1}((\tau,T],V)\text{, }u\in C^{0,\frac{1}{2}}([\tau+\varepsilon
,T],H^{2}(\Omega),\ \frac{d^{2}u}{dt^{2}}\in L^{2}(\tau+\varepsilon,T;H)\text{
for all }\varepsilon>0. \label{RegularityU}%
\end{equation}

\textbf{Step 2. }Let $u_{\tau}\in H.$

Since $u\in L^{2}(\tau+\varepsilon,T;D(A))$ for any $\varepsilon>0$, for any
$\delta>0$ there is $t_{0}\in\left(  \tau,\tau+\delta\right)  $ such that
$u\left(  t_{0}\right)  \in D(A)$. We know by Step 1 that $u\in C^{1}%
((t_{0},T],V)$ and $u\in C^{0,\frac{1}{2}}([t_{0}+\varepsilon,T],H^{2}%
(0,1)),\ \dfrac{d^{2}u}{dt^{2}}\in L^{2}(t_{0}+\varepsilon,T;H)$ for all
$\varepsilon>0$. As $\delta>0$ is arbitrary, we get (\ref{RegularityU}).

\textbf{Step 3. }Let $f\in C((\tau,T],H^{1}(0,1))$ and $u_{\tau}\in H.$

In such a case we have%
\[
Au=f-\frac{du}{dt}\in C((0,T],H^{1}(0,1)),
\]
so $u\in C((\tau,T],H^{3}(\Omega))$.

As $T>\tau$ is arbitrary, the result follows.
\end{proof}

\begin{corollary}
\label{Regularity2}If $f\in C((\tau,+\infty),H^{1}(0,1)),\ f\in L_{loc}%
^{2}\left(  \tau,+\infty;H\right)  \cap W_{loc}^{1,2}(\tau+\varepsilon
,+\infty;H)$ for all $\varepsilon>0$, then the partial derivatives
$u_{t},u_{xx}$ exists in the classical sense and are continuous on
$(\tau,+\infty)\times\lbrack0,1]$.
\end{corollary}

\begin{proof}
The continuous embedding $H^{1}\left(  0,1\right)  \subset C([0,1])$ implies
that%
\[
u_{xx},u_{t}\in C((\tau,+\infty),C([0,1]).
\]

\end{proof}

\bigskip

As it is known, increasing the temporal regularity of the function $f$ we can
prove that the solution $u$ is as regular as we desire.

\begin{lemma}
\label{Regularity3}Let $k\in\mathbb{N}$. Assume that $f\in L_{loc}^{2}\left(
\tau,+\infty;H\right)  \cap W_{loc}^{k+1,2}(\tau+\varepsilon,+\infty;H)$, for
all $\varepsilon>0$, and that%
\begin{equation}
f\in\cap_{j=0}^{k}C^{j}(\left(  \tau,+\infty),H^{2(k-j)+1}(0,1)\right)  .
\label{Condf}%
\end{equation}
Then, for any $u_{\tau}\in H$ the solution $u\left(  \text{\textperiodcentered
}\right)  $ to problem (\ref{Aux}) satisfies:
\begin{align}
u  &  \in C^{k+1}(\left(  \tau,+\infty\right)  ,V),\ u\in\cap_{j=0}^{k}%
C^{j}\left(  (\tau,+\infty),H^{2(k-j)+3}(0,1)\right)  ,\label{RegularityU2}\\
\frac{d^{k+2}u}{dt^{k+2}}  &  \in L_{loc}^{2}(\tau+\varepsilon,+\infty
;H)\text{ for all }\varepsilon>0.\nonumber
\end{align}

\end{lemma}

\begin{proof}
First, let us consider the case $k=1$, so $f\in L_{loc}^{2}\left(
\tau,+\infty;H\right)  \cap W_{loc}^{2,2}(\tau+\varepsilon,+\infty;H)$, for
all $\varepsilon>0$, and $f\in C^{1}(\left(  \tau,+\infty),H^{1}(0,1)\right)
\cap C(\left(  \tau,+\infty),H^{3}(0,1)\right)  $.

We know by Proposition \ref{Regularity} that $u\in C^{1}((\tau,+\infty),V)\cap
C((\tau,+\infty),H^{3}(0,1))$. The function $z_{1}=\dfrac{du}{dt}$ is the
unique solution to the problem%
\[
\left\{
\begin{array}
[c]{l}%
\dfrac{\partial z_{1}}{\partial t}-\dfrac{\partial^{2}z_{1}}{\partial x^{2}%
}=\dfrac{df}{dt},\text{ on }(\tau+\varepsilon,+\infty)\times(0,1),\\
z_{1}(t,0)=z_{1}(t,1)=0,\\
z_{1}(\tau+\varepsilon,x)=f(\tau+\varepsilon)-Au(\tau+\varepsilon)\in
H^{1}(0,1)\subset H,
\end{array}
\right.
\]
for any $\varepsilon>0$. Making use again of Proposition \ref{Regularity} and
taking into account that $\varepsilon$ is arbitrarily small we infer that
\begin{align*}
z_{1}  &  \in C^{1}((\tau,+\infty),V)\cap C((\tau,+\infty),H^{3}(0,1)),\\
\frac{d^{2}z_{1}}{d^{2}t}  &  \in L_{loc}^{2}(\tau+\varepsilon,+\infty
;H)\text{ for all }\varepsilon>0,
\end{align*}
so%
\begin{align*}
u  &  \in C^{2}((\tau,+\infty),V)\cap C^{1}((\tau,+\infty),H^{3}(0,1)),\\
\frac{d^{3}u}{dt^{3}}  &  \in L_{loc}^{2}(\tau+\varepsilon,+\infty;H)\text{
for all }\varepsilon>0.
\end{align*}
Finally,%
\[
Au=f-\frac{du}{dt}\in C(\left(  \tau,+\infty),H^{3}(0,1)\right)  ,
\]
so%
\[
u\in C(\left(  \tau,+\infty),H^{5}(0,1)\right)  .
\]

By induction, assume that the result is true for $k-1$. Hence, $u\in\cap
_{j=0}^{k-1}C^{j}(\left(  \tau,+\infty),H^{2(k-j)+1}(0,1)\right)  ,$ $u\in
C^{k}(\left(  \tau,+\infty\right)  ,V)$ and $\dfrac{d^{k+1}u}{dt^{k+1}}\in
L_{loc}^{2}(\tau+\varepsilon,+\infty;H)$ for all $\varepsilon>0$. The function
$z_{k}=\dfrac{d^{k}u}{dt^{k}}$ is the unique solution to the problem%
\[
\left\{
\begin{array}
[c]{l}%
\dfrac{\partial z_{k}}{\partial t}-\dfrac{\partial^{2}z_{k}}{\partial x^{2}%
}=\dfrac{d^{k}f}{dt^{k}},\text{ on }(\tau+\varepsilon,+\infty)\times(0,1),\\
z_{k}(t,0)=z_{k}(t,1)=0,\\
z_{k}(\tau+\varepsilon,x)=\dfrac{d^{k-1}f}{dt^{k-1}}(\tau+\varepsilon
)-A\dfrac{d^{k-1}u}{dt^{k-1}}(\tau+\varepsilon)\in H^{1}(0,1)\subset H,
\end{array}
\right.
\]
for any $\varepsilon>0$. Thus, Proposition \ref{Regularity} implies that%
\begin{align*}
z_{k}  &  \in C^{1}((\tau,+\infty),V)\cap C((\tau,+\infty),H^{3}(0,1)),\\
\frac{d^{2}z_{k}}{d^{2}t}  &  \in L_{loc}^{2}(\tau+\varepsilon,+\infty
;H)\text{ for all }\varepsilon>0,
\end{align*}
and then%
\begin{align*}
u  &  \in C^{k+1}((\tau,+\infty),V)\cap C^{k}((\tau,+\infty),H^{3}(0,1)),\\
\frac{d^{k+2}u}{dt^{k+2}}  &  \in L_{loc}^{2}(\tau+\varepsilon,+\infty
;H)\text{ for all }\varepsilon>0.
\end{align*}

Finally, we shall prove that $u\in C^{j}\left(  (\tau,+\infty),H^{2(k-j)+3}%
(0,1)\right)  $, for all $j\in\{0,1,...,k\}$ by induction. Assume that it is
true for $j\leq i\leq k$, where $j\in\{1,...,k\}$. Then%
\[
Az_{j-1}=\dfrac{d^{j-1}f}{dt^{j-1}}-\frac{dz_{j-1}}{dt}\in C\left(
(\tau,+\infty),H^{2(k-j)+3}(0,1)\right)  ,
\]
so%
\[
u\in C^{j-1}\left(  (\tau,+\infty),H^{2(k-(j-1))+3}(0,1)\right)  .
\]

\end{proof}

\begin{corollary}
\label{Regularity4}Under the conditions of Lemma \ref{Regularity3}, the
solution $u$ belongs to $C^{k+1}(\left(  \tau,+\infty\right)  \times
\lbrack0,1]).$
\end{corollary}

\begin{proof}
Since $\dfrac{d^{j}u}{dt^{j}}\in C\left(  (\tau,+\infty),H^{2(k-j)+3}%
(0,1)\right)  $, for all $j\in\{0,1,...,k+1\}$, we obtain that%
\[
\frac{\partial^{j+2\left(  k-j\right)  +2}u}{\partial t^{j}\partial
x^{2(k-j)+2}}\in C(\left(  \tau,+\infty\right)  ,C([0,1]).
\]
As $-j+2k+2\geq k+1$ for any $j\in\{0,...,k+1\}$, we infer that all the
partial derivatives of order less or equal to $k+1$ are continuous in $\left(
\tau,+\infty\right)  \times\lbrack0,1]$.
\end{proof}

\begin{corollary}
\label{Regularity5}If $f\in C^{\infty}\left(  \left(  \tau,+\infty\right)
\times\lbrack0,1]\right)  $ and $f\in L_{loc}^{2}\left(  \tau,+\infty
;H\right)  $, then the solution $u$ to problem (\ref{Aux}) belongs to
$C^{\infty}\left(  \left(  \tau,+\infty\right)  \times\lbrack0,1]\right)  $ as well.
\end{corollary}

\bigskip

\textbf{Acknowledgments.}

The author was partially supported by the Spanish Ministry of Science,
Innovation and Universities, project PGC2018-096540-B-I00, by the Spanish
Ministry of Science and Innovation, project PID2019-108654GB-I00, and by Junta
de Andaluc\'{\i}a (Spain), projects P18-FR-4509 and P18-FR-2025.

We would like to thank the anonymous referee for his/her useful remarks and comments.

\end{document}